\global \def \auxfile {\FALSE }
  \input amssym
  \input miniltx
  \input pictex

  \font \bbfive = bbm5
  \font \bbseven = bbm7
  \font \bbten = bbm10
  \font \eightbf = cmbx8
  \font \eighti = cmmi8 \skewchar \eighti = '177
  \font \eightit = cmti8
  \font \eightrm = cmr8
  \font \eightsl = cmsl8
  \font \eightsy = cmsy8 \skewchar \eightsy = '60
  \font \eighttt = cmtt8 \hyphenchar \eighttt = -1

  \font \sixi = cmmi6 \skewchar \sixi = '177
  \font \sixrm = cmr6
  \font \sixsy = cmsy6 \skewchar \sixsy = '60
  \font \tensc = cmcsc10

  \scriptfont \bffam = \bbseven
  \scriptscriptfont \bffam = \bbfive
  \textfont \bffam = \bbten

  \newskip \ttglue

  \def \eightpoint {\def \rm {\fam 0 \eightrm }\relax
  \textfont 0= \eightrm
  \scriptfont 0 = \sixrm \scriptscriptfont 0 = \fiverm
  \textfont 1 = \eighti
  \scriptfont 1 = \sixi \scriptscriptfont 1 = \fivei
  \textfont 2 = \eightsy
  \scriptfont 2 = \sixsy \scriptscriptfont 2 = \fivesy
  \textfont 3 = \tenex
  \scriptfont 3 = \tenex \scriptscriptfont 3 = \tenex
  \def \it {\fam \itfam \eightit }\relax
  \textfont \itfam = \eightit
  \def \sl {\fam \slfam \eightsl }\relax
  \textfont \slfam = \eightsl
  \def \bf {\fam \bffam \eightbf }\relax
  \textfont \bffam = \bbseven
  \scriptfont \bffam = \bbfive
  \scriptscriptfont \bffam = \bbfive
  \def \tt {\fam \ttfam \eighttt }\relax
  \textfont \ttfam = \eighttt
  \tt \ttglue = .5em plus.25em minus.15em
  \normalbaselineskip = 9pt
  \def \MF {{\manual opqr}\-{\manual stuq}}\relax
  \let \sc = \sixrm
  \let \big = \eightbig
  \setbox \strutbox = \hbox {\vrule height7pt depth2pt width0pt}\relax
  \normalbaselines \rm }

  \def \setfont #1{\font \auxfont =#1 \auxfont }
  \def \withfont #1#2{{\setfont {#1}#2}}

  \def \TRUE {Y}
  \def \FALSE {N}
  \def \EMPTY {}

  \def \ifundef #1{\expandafter \ifx \csname #1\endcsname \relax }

  \def \undefrule {\kern 2pt \vrule width 5pt height 5pt depth 0pt \kern 2pt}
  \def \UndefLabels {}
  \def \possundef #1{\ifundef {#1}\undefrule {\eighttt #1}\undefrule
    \global \edef \UndefLabels {\UndefLabels #1\par }
  \else \csname #1\endcsname \fi }

  \newcount \secno \secno = 0
  \newcount \stno \stno = 0
  \newcount \eqcntr \eqcntr = 0

  \ifundef {showlabel} \global \def \showlabel {\FALSE } \fi
  \ifundef {auxfile} \global \def \auxfile {\TRUE } \fi

  \def \define #1#2{\global \expandafter \edef \csname #1\endcsname {#2}}
  \long \def \error #1{\medskip \noindent {\bf ******* #1}}
  \def \fatal #1{\error {#1\par Exiting...}\end }

  \def \advseqnumbering {\global \advance \stno by 1 \global \eqcntr =0}

  \def \current {\ifnum \secno = 0 \number \stno \else \number \secno \ifnum
\stno = 0 \else .\number \stno \fi \fi }

  \begingroup \catcode `\@ =0 \catcode `\\=11 @global @def @textbackslash {\}
@endgroup
  \def \space { }

  \def \deflabel #1#2{\if \TRUE \showlabel \hbox {\sixrm [[ #1 ]]} \fi
    \ifundef {#1PrimarilyDefined}\define {#1}{#2}\define
{#1PrimarilyDefined}{#2}\if \TRUE \auxfile \immediate \write 1 {\textbackslash
newlabel {#1}{#2}}\fi
    \else
      \edef \old {\csname #1\endcsname }\edef \new {#2}\if \old \new \else
\fatal {Duplicate definition for label ``{\tt #1}'', already defined as ``{\tt
\old }''.}\fi
      \fi }

  \def \newlabel #1#2{\define {#1}{#2}}
  \def \label #1 {\deflabel {#1}{\current }}

  \def \equationmark #1 {\ifundef {InsideBlock}
	  \advseqnumbering
	  \eqno {(\current )}
	  \deflabel {#1}{\current }
	\else
	  \global \advance \eqcntr by 1
	  \edef \subeqmarkaux {\current .\number \eqcntr }
	  \eqno {(\subeqmarkaux )}
	  \deflabel {#1}{\subeqmarkaux }
	\fi }

  \def \split #1.#2.#3.#4;{\global \def \parone {#1}\global \def \partwo
{#2}\global \def \parthree {#3}\global \def \parfour {#4}}
  \def \NA {NA}
  \def \ref #1{\split #1.NA.NA.NA;(\possundef {\parone }\ifx \partwo \NA \else
.\partwo \fi )}

  \newcount \bibno \bibno = 0
  \def \newbib #1#2{\define {#1}{#2}}

  \def \Bibitem #1 #2; #3; #4 \par {\smallbreak
    \global \advance \bibno by 1
    \item {[\possundef {#1}]} #2, {``#3''}, #4.\par
    \ifundef {#1PrimarilyDefined}\else
      \fatal {Duplicate definition for bibliography item ``{\tt #1}'', already
defined in ``{\tt [\csname #1\endcsname ]}''.}
      \fi
	\ifundef {#1}\else
	  \edef \prevNum {\csname #1\endcsname }
	  \ifnum \bibno =\prevNum \else
		\error {Mismatch bibliography item ``{\tt #1}'', defined earlier (in aux file
?) as ``{\tt \prevNum }'' but should be
	``{\tt \number \bibno }''.  Running again should fix this.}
		\fi
	  \fi
    \define {#1PrimarilyDefined}{#2}\if \TRUE \auxfile \immediate \write 1
{\textbackslash newbib {#1}{\number \bibno }}\fi }

  \def \jrn #1, #2 (#3), #4-#5;{{\sl #1}, {\bf #2} (#3), #4--#5}
  \def \Article #1 #2; #3; #4 \par {\Bibitem #1 #2; #3; \jrn #4; \par }

  \def \references {\begingroup \bigbreak \eightpoint \centerline {\tensc
References} \nobreak \medskip \frenchspacing }

  \catcode `\@ =11
  \def \c@itrk #1{{\bf \possundef {#1}}}
  \def \c@ite #1{{\rm [\c@itrk {#1}]}}
  \def \sc@ite [#1]#2{{\rm [\c@itrk {#2}\hskip 0.7pt:\hskip 2pt #1]}}
  \def \du@lcite {\if \pe@k [\expandafter \sc@ite \else \expandafter \c@ite \fi
}
  \def \cite {\futurelet \pe@k \du@lcite }
  \catcode `\@ =12

  \def \Headlines #1#2{\nopagenumbers
    \headline {\ifnum \pageno = 1 \hfil
    \else \ifodd \pageno \tensc \hfil \lcase {#1} \hfil \folio
    \else \tensc \folio \hfil \lcase {#2} \hfil
    \fi \fi }}

  \def \title #1{\medskip \centerline {\withfont {cmbx12}{\ucase {#1}}}}

  \long \def \Quote #1\endQuote {\begingroup \leftskip 35pt \rightskip 35pt
\parindent 17pt \eightpoint #1\par \endgroup }
  \long \def \Abstract #1\endAbstract {\vskip 1cm \Quote \noindent #1\endQuote }
  
  \def \Authors #1{\bigskip \centerline {\tensc #1}}
  \def \Note #1{\footnote {}{\eightpoint #1}}
  \def \Date #1 {\Note {\it Date: #1.}}

  \def \today {\ifcase \month \or January\or February\or March\or April\or
May\or June\or July\or August\or September\or
    October\or November\or December\fi \space \number \day , \number \year }
  \def \hoje {\number \day /\number \month /\number \year }

  \newcount \auxone \newcount \auxtwo \newcount \auxthree
  \def \fulltime {\auxone =\time \auxtwo =\time \divide \auxone by 60 \auxthree
=\auxone \multiply \auxthree by 60 \advance
    \auxtwo by -\auxthree \hoje \ \ \ifnum \auxone <10 0\fi \number \auxone
:\ifnum \auxtwo <10 0\fi \number \auxtwo }

  \def \part #1#2{\vfill \eject \null \vskip 0.3\vsize
    \withfont {cmbx10 scaled 1440}{\centerline {PART #1} \vskip 1.5cm
\centerline {#2}} \vfill \eject }

  \def \fix {\smallskip \noindent $\blacktriangleright $\kern 12pt}
  \def \iskip {\medskip \noindent }

  \def \ucase #1{\edef \auxvar {\uppercase {#1}}\auxvar }
  \def \lcase #1{\edef \auxvar {\lowercase {#1}}\auxvar }
  \def \section #1 \par {\global \advance \secno by 1 \stno = 0
    \goodbreak \bigbreak
    \noindent {\bf \number \secno .\enspace #1.}
    \nobreak \medskip \noindent }
  \def \state #1 #2\par {\begingroup \def \InsideBlock {} \medbreak \noindent
\advseqnumbering {\bf \current .\enspace  #1.\enspace \sl #2\par }\medbreak
\endgroup }

  \def \definition #1\par {\state Definition \rm #1\par }

  \newcount \CloseProofFlag
  \def \closeProof {\eqno \endproofmarker \global \CloseProofFlag =1}
  \long \def \Proof #1\endProof {\begingroup \def \InsideBlock {} \global
\CloseProofFlag =0
     \medbreak \noindent {\it Proof.\enspace }#1
     \ifnum \CloseProofFlag =0 \hfill $\endproofmarker $ \looseness = -1 \fi
     \medbreak \endgroup }

  \def \quebra #1{#1 $$$$ #1}
  \def \explica #1#2{\mathrel {\buildrel \hbox {\sixrm #1} \over #2}}
  \def \explain #1#2{\explica {\ref {#1}}{#2}}
  
  \def \=#1{\explain {#1}{=}}

  \def \pilar #1{\vrule height #1 width 0pt}
  \def \stake #1{\vrule depth #1 width 0pt}

  \newcount \fnctr \fnctr = 0
  \def \fn #1{\global \advance \fnctr by 1
    \edef \footnumb {$^{\number \fnctr }$}\footnote {\footnumb }{\eightpoint
#1\par \vskip -10pt}}

  \def \text #1{\hbox {#1}}
  \def \bool #1{[{\scriptstyle #1}]\,}

  \def \Item #1{\smallskip \item {{\rm #1}}}
  \newcount \zitemno \zitemno = 0

  \def \izitem {\global \zitemno = 0}
  \def \zitemplus {\global \advance \zitemno by 1 \relax }
  \def \rzitem {\romannumeral \zitemno }
  \def \rzitemplus {\zitemplus \rzitem }
  \def \zitem {\Item {{\rm (\rzitemplus )}}}

  \def \zitemmark #1 {\deflabel {#1}{\current .\rzitem } \deflabel
{Local#1}{\rzitem }}

  \newcount \nitemno \nitemno = 0
  
  \def \nitem {\global \advance \nitemno by 1 \Item {{\rm (\number \nitemno )}}}

  \newcount \aitemno \aitemno = -1
  \def \boxlet #1{\hbox to 6.5pt{\hfill #1\hfill }}
  \def \iaitem {\aitemno = -1}
  \def \aitemconv {\ifcase \aitemno a\or b\or c\or d\or e\or f\or g\or  h\or
i\or j\or k\or l\or m\or n\or o\or p\or q\or r\or s\or t\or u\or  v\or w\or x\or
y\or z\else zzz\fi }
  \def \aitem {\global \advance \aitemno by 1\Item {(\boxlet \aitemconv )}}
  \def \aitemmark #1 {\deflabel {#1}{\aitemconv }}

  \def \Case #1:{\medskip \noindent {\tensc Case #1:}}

  \def \<{\left \langle \vrule width 0pt depth 0pt height 8pt }
  \def \>{\right \rangle }
  \def \({\big (}
  \def \){\big )}
  \def \ds {\displaystyle }
  \def \and {\hbox {\quad and \quad }}

  \def \imply {\mathrel {\Rightarrow }}
  \def \IFF {\kern 7pt\Leftrightarrow \kern 7pt}
  \def \IMPLY {\kern 7pt \Rightarrow \kern 7pt}
  \def \for #1{\quad \forall \,#1}
  \def \endproofmarker {\square }
  \def \"#1{{\it #1}\/} \def \umlaut #1{{\accent "7F #1}}
  \def \inv {^{-1}}
  \def \*{\otimes }
  \def \caldef #1{\global \expandafter \edef \csname #1\endcsname {{\cal #1}}}
  \def \bfdef #1{\global \expandafter \edef \csname #1\endcsname {{\bf #1}}}
  \bfdef N \bfdef Z \bfdef C \bfdef R

  \if \TRUE \auxfile
    \IfFileExists {\jobname .aux}{\input \jobname .aux}{\null }
    \immediate \openout 1 \jobname .aux
    \fi

  \def \close {\if \EMPTY \UndefLabels \else
      \message {*** There were undefined labels ***} \iskip
      ****************** \ Undefined Labels: \tt \par \UndefLabels
      \fi
    \if \TRUE \auxfile \closeout 1 \fi
    \par \vfill \supereject \end }

  \def \inpatex #1.tex{\input #1.atex}

  \def \Caixa #1{\setbox 1=\hbox {$#1$\kern 1pt}\global \edef \tamcaixa {\the
\wd 1}\box 1}
  \def \caixa #1{\hbox to \tamcaixa {$#1$\hfil }}

  \def \med #1{\mathop {\textstyle #1}\limits }
  
  \def \medsum {\med \sum }
  
  \def \medcup {\med \bigcup }
  \def \medcap {\med \bigcap }

   \newlabel {StandingHyp}{2.1} \newlabel {FormulaProduct}{2.2} \newlabel
{FormulaProductSimpler}{2.3} \newlabel {InductionSect}{3} \newlabel
{BasicNotations}{3.1} \newlabel {DefineForm}{3.2} \newlabel
{InnProdRightLinear}{3.3} \newlabel {QBasis}{3.4} \newlabel
{DefinEinducedModule}{3.6} \newlabel {IndPreserveIrred}{3.7} \newlabel
{DeltaDeltaInMod}{3.7.1} \newlabel {VectorAnnnihiletor}{3.8} \newlabel
{IndIdealDependsOnlyIdeal}{3.9} \newlabel {IntroJXzeroI}{3.10} \newlabel
{AnnInducedVsIndIdeal}{3.11} \newlabel {DeltaKbDeltaL}{3.13} \newlabel
{UseDeltaKbDeltaL}{3.14} \newlabel {Examples}{3.15} \newlabel
{InterCoefAlg}{3.16} \newlabel {IntroE}{3.17} \newlabel {CondExp}{3.18}
\newlabel {IntroNu}{3.19} \newlabel {IntroF}{3.20} \newlabel {FvsBilForm}{3.21}
\newlabel {AltDescription}{3.22} \newlabel {AdmIdSect}{4} \newlabel
{FAlmostMult}{4.1} \newlabel {FunctOut}{4.2} \newlabel {FunctOutTwo}{4.3}
\newlabel {MapIdeals}{4.4} \newlabel {IntroIPrime}{4.5} \newlabel
{DefAdmiss}{4.6} \newlabel {ChangeIdeals}{4.7} \newlabel {TrivialIdeals}{4.8}
\newlabel {LiftC}{4.9} \newlabel {AdmissibleCharact}{4.10} \newlabel
{LargestIdeal}{4.11} \newlabel {FJIsAdmissible}{4.12} \newlabel {RepSec}{5}
\newlabel {IdealIsKernel}{5.1} \newlabel {LocalUnits}{5.2} \newlabel
{BOneC}{5.2.1} \newlabel {LcNonDeg}{5.3} \newlabel {Disintegration}{5.4}
\newlabel {CovarCond}{5.4.i} \newlabel {LocalCovarCond}{i} \newlabel
{Reconstruct}{5.4.ii} \newlabel {LocalReconstruct}{ii} \newlabel
{Projecao}{5.4.iii} \newlabel {LocalProjecao}{iii} \newlabel {DefineUg}{5.4.1}
\newlabel {EigenVec}{5.5} \newlabel {qxpifphiv}{5.5.1} \newlabel
{SummarizePix}{5.6} \newlabel {MuOnEg}{5.8} \newlabel {DefineBigUg}{5.9}
\newlabel {Functoriality}{5.10} \newlabel {Discretization}{5.12} \newlabel
{IntroQ}{5.13} \newlabel {LocalValueZero}{5.13.1} \newlabel {XiSuppD}{5.13.2}
\newlabel {FirstNulSpace}{5.14} \newlabel {SecondNulSpace}{5.15} \newlabel
{FormulasOnQx}{5.16} \newlabel {ExpressionForPixU}{5.17} \newlabel
{MatrixEntries}{5.18} \newlabel {SoFarKernelInclusions}{5.19} \newlabel
{ThirdNulSpace}{5.20} \newlabel {WellDef}{5.20.1} \newlabel
{HypothesisPiVanishes}{5.20.2} \newlabel {MeanwhileClaim}{5.20.3} \newlabel
{OrbitInvariant}{5.21} \newlabel {DefineRepOrbit}{5.22} \newlabel
{JIsIntersRho}{5.23} \newlabel {RhoSect}{6} \newlabel {EnterKhModule}{6.1}
\newlabel {RhoIsInduced}{6.2} \newlabel {ThisIsCovariance}{6.2.1} \newlabel
{RhofDgOnTau}{6.2.2} \newlabel {MainResult}{6.3} \newlabel {FIntersection}{7.2}
\newlabel {TopFreeSect}{8} \newlabel {DefTopFree}{8.1} \newlabel
{TopFreeOnOrb}{8.2} \newlabel {StrongerProperty}{8.3} \newlabel
{TopFreeNoAdmiss}{8.4} \newlabel {FixingConjuga}{9.3} \newlabel
{LemaOnMembership}{9.4} \newlabel {LemaOnMembershipTwo}{9.4.1} \newlabel
{ConditionForRegAdmiss}{9.5} \newlabel {CommutativeAdmiss}{9.6} \newlabel
{NormalIsAdmissible}{10.3} \newlabel {SourceOfNormal}{10.4} \newlabel
{AugmAdmiss}{10.5} \newlabel {TranspositionSect}{11} \newlabel
{DefineTransp}{11.1} \newlabel {InclusionInduced}{11.2} \newlabel
{MainEqualityInduced}{11.3} \newlabel {DescribeTransposition}{11.4} \newlabel
{CondOnV}{11.4.1} \newlabel {OneVatATime}{11.4.2} \newbib {LisaEtAll}{1} \newbib
{ReznikofEtAll}{2} \newbib {EH}{3} \newbib {mybook}{4} \newbib {Dan}{5} \newbib
{GR}{6} \newbib {IonWill}{7} \newbib {Renault}{8} \newbib {Sauvageot}{9} \newbib
{FS}{10} \newbib {SimsWill}{11} \newbib {SteinbAlg}{12} \newbib {Steinberg}{13}
\newbib {Willard}{14} \newbib {WillBook}{15}

\font \rs = rsfs10  \def \curlyL {\hbox {\rs
L}} \def \lc #1{\curlyL _c(#1)} \def \Lc {\lc X} \def \cp #1{\Lc {\,\rtimes \,}
#1} \def \alg {\cp G} \def \M {M} \def \bool #1{\,[{\raise 0.5pt \hbox
{$\scriptstyle #1$}}]\,}

\def \indbox {\hbox {\sl Ind}} \def \ind #1#2{\indbox \def \a {}\def \b {#1}\ifx
\a \b \else _{#1}\fi (#2)} \def \JxComXZeroI {\ind {x_0}I} \def \Jx #1{\ind
{}{#1}} \def \JxBig #1{\indbox \Big (#1\Big )} \def \JxI {\Jx I}

\def \ideal {\mathrel \trianglelefteq }   \def \supp {\hbox {supp}} \def
\thickdot {\hbox {\withfont {cmbx12}{.}}} \def \d {\Delta } \def \Orb
{\mathchoice {\hbox {Orb}}{\hbox {Orb}}{\hbox {\eightrm Orb}}{\hbox {\sixrm
Orb}}} \def \Ker {\mathchoice {\hbox {Ker}}{\hbox {Ker}}{\hbox {\eightrm
Ker}}{\hbox {\sixrm Ker}}} \def \cOrb {\overline {\pilar {8pt}\Orb }} \def \Sx
{{\cal S}} \def \Hx {H}  \def \absform {\langle \cdot \kern
1.3pt,\cdot \rangle } \def \lin {L}

\def \Fld {K} \def \GpAlg #1{\Fld \!#1} \def \KH {\GpAlg {\Hx }} \def \subsub
#1#2{_{\buildrel \scriptstyle #1 \over {\pilar {7pt}#2}}} \def \abar #1#2{\bar
\alpha _{#1}(#2)}  \def \Flecha #1#2{\ \buildrel \textstyle #2 \over {\hbox to
#1 {\rightarrowfill }}\ } \def \flecha #1{\setbox 0 \hbox {$\textstyle #1$}
\dimen 0=\wd 0 \advance \dimen 0 by 20pt \Flecha {\dimen 0}{#1}} \def \medcap
#1{\raise 2pt \hbox {$\mathop {\textstyle \bigcap }\limits _{#1}$}} \def \klopen
{compact-open} \def \book #1#2{\cite [#1 #2]{mybook}}

\Headlines {The ideal structure of algebraic partial crossed products}{M. Dokuchaev and R. Exel} \title {The ideal
structure of algebraic} \title {partial crossed products}

\Authors {M. Dokuchaev and R. Exel} \Date {22 May 2014}

\Note {\it Key words and phrases: \rm Effros-Hahn conjecture, partial crossed
product, induced representation, induced ideal, locally constant function.}
\Note {Partially supported by CNPq and FAPESP.}

\bigskip  \Abstract
  \def \Lc {\hbox {\withfont {rsfs10 scaled 833}L}_c(X)}Given a partial action
of a discrete group $G$ on a Hausdorff, locally compact, totally disconnected
topological space $X$, we consider the correponding partial action of $G$ on the
algebra $\Lc $ consisting of all locally constant, compactly supported functions
on $X$, taking values in a given field $\Fld $.  We then study the ideal
structure of the algebraic partial crossed product $\alg $.  After developping a
theory of induced ideals, we show that every ideal in $\alg $ may be obtained as
the intersection of ideals induced from isotropy groups, thus proving an
algebraic version of the Effros-Hahn conjecture.
  \endAbstract

\section Introduction

The study of ideals in crossed product C*-algebras has a long history and is
best subsumed by the quest to prove and generalize the celebrated Effros-Hahn
conjecture \cite {EH} formulated roughly fifty years ago.  In its original form,
the conjecture states that every primitive ideal in the crossed product of a
commutative C*-algebra by a locally compact group should be induced from a
primitive ideal in the C*-algebra of some isotropy group.

For the case of discrete amenable groups, the Effros-Hahn conjecture was proven
by Sauvageot \cite {Sauvageot}, and since then has been extended to various
other contexts, notably to locally compact groups acting on non commutative
C*-algebras, as proven by Gootman and Rosenberg \cite {GR} under separability
conditions.
  Motivated by Fack and Skandalis' study of C*-algebras associated to foliations
\cite {FS}, Renault \cite {Renault} realized that the Effros-Hahn conjecture
also applies for groupoid C*-algebras and proved a version of it in this
context.  This was later refined by Ionescu and Williams in \cite {IonWill}.

Most treatments of the Effros-Hahn conjecture focus on the conjecture itself,
namely describing a previously given ideal in the crossed product algebra in
terms of induced ideals, rather than studiyng the relationship between the
\"{input} ideal in the isotropy group algebra and its corresponding \"{output}
induced ideal.

To be fair, there are a few works in the literature where this delicate
relationship is discussed.  Among them we should mention \cite [Theorem
8.39]{WillBook}, where a complete classification is given for the collection of
primitive ideals on $C_0(X)\rtimes G$, where $X$ is locally compact and $G$
abelian.  There, it is shown that every such ideal is induced from some isotropy
group and hence arises from a pair $(x,\chi )$, where $x$ is a point in $X$, and
$\chi $ a character on $G$. Most importantly, the primitive ideals for two pairs
$(x_1,\chi _1)$ and $(x_2,\chi _2)$ coincide if and only if the closure of the
orbit of $x_1$ coincides with that of $x_2$, and $\chi _1$ coincides with $\chi
_2$ on the isotropy group of $x_1$.

In the above situation, when two points have the same orbit closure, it may be
shown that their isotropy groups coincide, but this relies heavily on the
commutativity of $G$.  In particular it is not even clear how to phrase the
above condition in case $G$ is not commutative.

The question of when two induced ideals coincide is also taken up in \cite
{SimsWill}, where the object of study is the C*-algebra of the Deaconu-Renault
groupoid built from an action of the semigroup ${\bf N}^k$ by local
homeomorphisms on a locally compact space $X$.  The primitive ideals of this
C*-algebra are shown to be parametrized by pairs $(x,\chi )$, where $x$ is a
point in $X$, and $\chi $ a character on ${\bf N}^k$, i.e.~an element of ${\bf
T}^k$, and again a criterion is given for when two such pairs
  lead to the same primitive ideal.  The condition is that the two points of $X$
must have identical orbit closures and the two characters must coincide as
characters on the interior of the isotropy of the groupoid reduced to the common
orbit closure.  In a sense, this result also relies on commutativity.

In the present paper our aim is to study the Effros-Hahn conjecture in a new
setting, namely the algebraic partial crossed product $\alg $, where $G$ is a
not necessarily commutative, discrete group partially acting on a locally
compact, totally disconnected topological space $X$, and $\Lc $ is the algebra
consisting of all locally constant, compactly supported functions on $X$, taking
values in a given field $\Fld $.

The justification for studying this setting comes from the current interest to
investigate purely algebraic versions os some intensely studied C*-algebras,
such as the Leavitt path algebras which may be viewed as algebraic counterparts
of graph C*-algebras.  In many such cases the pertinent C*-algebra is a
C*-algebraic partial crossed product of the form $C_0(X)\rtimes G$, with $X$
totally disconnected, while its algebraic sibling is the algebraic partial
crossed product $\alg $.  Steinberg algebras \cite {SteinbAlg} in fact
generalize this correspondence to the case where an ample \'etale groupoid
replaces the above partial action.

One of the main results of the present paper, namely Theorem \ref {MainResult},
is a version of the Effros-Hahn conjecture, where we prove that every ideal of
$\alg $ is given as the intersection of ideals induced from isotropy groups. The
method of proof is entirely elementary and does not rely on the measure
theoretical or analytical tools on which the main proofs in \cite {Sauvageot},
\cite {GR} and \cite {Renault} are based, chiefly because our setting is
eminently algebraic. The strategy adopted here is as follows:
  given an ideal $J$ of $\alg $, we first choose a representation $\pi $ of
$\alg $ whose null space coincides with $J$.  We then build another
representation, which we call the \"{discretization} of $\pi $,
  whose null space coincides with that of $\pi $, and hence also with $J$.
  The discretized representation is then easily seen to decompose as a direct
sum of sub-representations based on the orbits for the action of $G$ on $X$.
Each such sub-representation is finally shown to be equivalent to an induced
representation, and hence the initially given ideal $J$ is seen to coincide with
the intersection of the null spaces of the various induced representation
involved, each of which is then an induced ideal.

Only a tiny amount of the theory of induced ideals is necessary to prove Theorem
\ref {MainResult}, our version of the Effros-Hahn conjecture, but still we have
chosen to start the study of induced ideals before stating and proving \ref
{MainResult}, mostly in order to be able to refer to the main concepts involved.

Before and after the proof of Theorem \ref {MainResult}, in fact throughout the
paper, we develop tools designed to understand the induction process itself,
attempting to describe how exactly does an induced ideal $\JxComXZeroI $ depends
on point $x_0$ and on the ideal $I$ it is induced from.  From the outset, this
dependency is expected to be quite tricky for the following reason: there are
numerous examples where a crossed product algebra turns out to be simple (see
e.g.~\cite [Theorem 4.1]{Dan}), and hence the assortment of ideals in the
crossed product is rather boring, but still there may be points with nontrivial
isotropy (not too many since the action must be topologically free by \cite
{LisaEtAll}, but this does not ruled out all points), and hence there may be
many ideals presenting themselves as input for the induction process. However,
as already mentioned, the output could be totally uninteresting due to
simplicity.

The explanation for this phenomenon is that, when inducing from the isotropy
group $\Hx _{x_0}$, where $x_0$ is some point in $X$, not all ideals in $\GpAlg
\Hx _{x_0}$ play a relevant role.  Those which do, namely the ones we call
\"{admissible}, are the only ones deserving attention in the sense that for
every ideal $I\ideal \GpAlg \Hx _{x_0}$, there exists a unique admissible ideal
$I'\subseteq I$, which induces the same ideal of $\alg $ as $I$ does.  This is
the content of Corollary \ref {ChangeIdeals}.
  The correspondence $I\mapsto \JxComXZeroI $ is consequently seen to be a
one-to-one mapping from the set of admissible ideals in $\GpAlg \Hx _{x_0}$ to
the set of ideals in $\alg $.

This naturally raises the question of which are the admissible ideals in $\GpAlg
\Hx _{x_0}$, a question we answer in general in \ref {AdmissibleCharact}, and
then in a few special cases in
  \ref {TopFreeNoAdmiss},
  \ref {ConditionForRegAdmiss},
  \ref {CommutativeAdmiss} and
  \ref {NormalIsAdmissible}.

We then consider the question of when two induced ideals $\ind {x_0}I$ and $\ind
{\tilde x_0}{\tilde I}$ coincide, where $x_0$ and $\tilde x_0$ are two points in
$X$, $I$ is an admissible ideal in $\GpAlg \Hx _{x_0}$, and $\tilde I$ is an
admissible ideal in $\GpAlg \Hx _{\tilde x_0}$ (based on our study of induced
ideals it suffices to consider the case where $I$ and $\tilde I$ are
admissible).

As already mentioned, when $x_0=\tilde x_0$ one has that
  $$
  \ind {x_0}I = \ind {\tilde x_0}{\tilde I} \ \iff \ I=\tilde I,
  $$
  so the remaining case is when $x_0\neq \tilde x_0$.  As in \cite {WillBook}
and \cite {SimsWill}, a necessary condition for
  $\ind {x_0}I$
  and
  $\ind {\tilde x_0}{\tilde I}$
  to coincide
  is that the orbits of $x_0$ and $\tilde x_0$ have the same closure (as long as
$I$ and $\tilde I$ are proper ideals, according to \ref {InclusionInduced}) but
now, in the absence of commutativity, the isotropy groups of $x_0$ and $\tilde
x_0$ no longer need to be the same, so comparing the extra data
  within
  $\GpAlg \Hx _{x_0}$ and
  $\GpAlg \Hx _{\tilde x_0}$
  (namely $I$ versus $\tilde I$)
  is no longer a straightforward matter.

To deal with this situation we introduce the notion of \"{transposition} of
ideals in \ref {DefineTransp} which is a way of comparing ideals in different
isotropy groups. Our main result in that direction, namely Theorem \ref
{MainEqualityInduced}, says that
  $
  \ind {x_0}I = \ind {\tilde x_0}{\tilde I}
  $
  if and only if $\tilde I$ is the transposition of $I$ and vice versa.

We have already mentioned that the algebra $\alg $, which is our main object of
study, may also be described as the Steinberg algebra for the transformation
groupoid associated to the partial action of $G$ on $X$.  Steinberg's results
obtained in \cite {SteinbAlg} and \cite {Steinberg} therefore apply to our
situation as well. On the other hand, in all likelihood our results may be shown
to hold for Steinberg algebras with minor modifications in our proofs.

Our algebras are all taken to be over a fixed field $\Fld $, but in most places
our results hold under the more general assumption that $\Fld $ is just a unital
commutative ring.  Notable exceptions are \ref {InterCoefAlg} and \ref
{TopFreeNoAdmiss}, where invertibility of nonzero elements in $\Fld $ is
crucial.

The second named author would like to acknowledge financial support from the
Funda\c c\~ao de Amparo \`a Pesquisa do Estado de S\~ao Paulo (FAPESP) during a
visit to the University of S\~ao Paulo, where a large part of this work was
conducted.  He would also like to acknowledge the warm hospitality of the
members of the Mathematics Department during that visit.

\section Preliminaries

Throughout most of this work we will assume the following:

\state {Standing Hypotheses} \label StandingHyp \rm
  \iaitem
  \aitem $\Fld $ is a field,
  \aitem $G$ is a discrete group,
  \aitem $X$ is a Hausdorff, locally compact, totally
  disconnected\fn {A locally compact topological space is
  \"{totally disconnected} if and only if it admits a basis of open sets
consisting of sets which are also closed \cite [Theorem 29.7]{Willard}.}
  topological space,
  \aitem $\stake {8pt}\theta =(\{\theta _g\}_{g\in G},\{X_g\}_{g\in G})$ is a
(topological) partial action \book {Definition}{5.1} of $G$ on $X$, such that
$X_g$ is \"{clopen} (closed and open) for every $g$ in $G$,
  \aitem whenever appropriate, we will also fix a distinguished point $x_0$ in
$X$.

\bigskip  Recall that a function
  $$
  f:X\to \Fld
  $$
  is said to be \"{locally constant} if, for every $x$ in $X$, there exists a
neighborhood $V$ of $x$, such that $f$ is constant on $V$.  The \"{support} of
$f$ is defined to be the set
  $$
  \supp (f) = \{x\in X: f(x)\neq 0\}.
  $$
  Observe that the support of a locally constant function $f$ is always closed,
so we will not bother to define the support as the \"{closure} of the above set,
as sometimes done in analysis.

By virtue of being locally compact and totally disconnected, we have that the
topology of $X$ admits a basis formed by
  \"{\klopen }\fn {That is, sets which are simultaneously compact and open.}
  subsets.  Given any {\klopen } set $E\subseteq X$, it is easy to see that its
characteristic function, here denoted by $1_E$, is locally constant and
compactly supported.  Moreover, one may easily prove that every locally
constant, compactly supported function $f:X\to \Fld $ is a linear combination of
the form
  $$
  f = \sum _{i=1}^nc_i\,1_{E_i},
  $$
  where the $E_i$ are pairwise disjoint {\klopen } subsets and the $c_i$ lie in
$\Fld $.

We will henceforth denote by $\Lc $ the set of all locally constant, compactly
supported, $\Fld $-valued functions on $X$.  With pointwise multiplication, $\Lc
$ is a commutative $\Fld $-algebra, which is unital if and only if $X$ is
compact.

For each $g$ in $G$, we may also consider the $\Fld $-algebra $\lc {X_g}$, which
we will identify with the set formed by all $f$ in $\Lc $ vanishing on
$X\setminus X_g$.  Under this identification $\lc {X_g}$ becomes an ideal in
$\Lc $.

Regarding the homemorphism $\theta _g:X_{g\inv }\to X_g$, we may define an
isomorphism
  $$
  \alpha _g:\lc {X_{g\inv }}\to \lc {X_g},
  $$
  by setting
  $$
  \alpha _g(f) = f\circ \theta _{g\inv }, \for f\in \lc {X_{g\inv }}.
  $$

  The collection formed by all ideals $\lc {X_g}$, together with the collection
of all $\alpha _g$, is then easily seen to be an (algebraic) partial action
\book {Definition}{6.4} of $G$ on $\Lc $.

This is a \"{unital} partial action (one for which the domain ideals are unital)
if and only if all of the $X_g$ are compact.  However we shall prefer to
consider the more general situation where the $X_g$ are only assumed to be
closed (besides being open).

The main goal of this paper is to study the algebraic crossed product
  $$
  \alg ,
  $$
  as defined in \book {Definition}{8.3}.  A general element $b\in \alg $ will be
denoted by
  $$
  b = \sum _{g\in G}f_g\d _g,
  $$
  where each $f_g$ lies in $\lc {X_{g\inv }}$, and $f_g=0$, for all but finitely
many group elements $g$.

In many texts dealing with crossed products the above place markers ``$\d _g$''
are denoted ``$\delta _g$'', but we shall reserve the latter to denote elements
in $\GpAlg G$, such as
  $$
  \sum _{h\in H}c_h\delta _h,
  $$
  where the $c_h$ are scalars in $\Fld $, again equal to zero except for
finitely many group elements $g$.  In fact, throughout this paper all summations
will be finite, either because the set of indices is finite, or because all but
finitely many summands are supposed to vanish.

\begingroup \def \a {e} \def \b {f} \def \g {g} \def \h {h}

Since we are asuming that every $X_g$ is clopen, its characteristic function
$1_{X_g}$, which we will abbreviate to
  $$
  1_g:=1_{X_g},
  $$
  is a locally constant function, although not necessarily compactly supported.
However, given any $f$ in $\Lc $, one has that $f1_{g\inv }$ is compactly
supported, so it belongs to $\lc {X_g}$.  We may therefore define
  $$
  \abar gf := \alpha _g(f1_{g\inv }).
  $$
  so that $\bar \alpha _g$ is a globally defined endomorphism of $\Lc $.

Recall that if $\g $ and $\h $ are elements of $G$, and if $\a \in \lc {X_\g }$,
and $\b \in \lc {X_\h }$, then the product of $\a \d _{\g }$ by $\b \d _{\h }$
is defined by
  $$
  (\a \d _{\g })(\b \d _{\h }) = \alpha _{\g }\big (\alpha _{\g \inv }(\a )\b
\big )\d _{\g \h }.
  \equationmark FormulaProduct
  $$

  In our present situation this expression may be made simpler as follows: since
$\alpha _{\g \inv }(\a )$ lies in $\lc  {X_{g\inv }}$, we have that
  $$
  \alpha _{\g \inv }(\a )=\alpha _{\g \inv }(\a ) 1_{\g \inv },
  $$
  so the coefficient of $\d _{\g \h }$ in \ref {FormulaProduct} equals
  $$
  \alpha _{\g }\big (\alpha _{\g \inv }(\a )\b \big ) =
  \alpha _{\g }\big (\alpha _{\g \inv }(\a )1_{\g \inv }\b \big ) =
  \alpha _{\g }\big (\alpha _{\g \inv }(\a )\big )\alpha _{\g }\big (1_{\g \inv
}\b \big ) =
  \a \abar \g \b .
  $$

  The promissed simpler formula for the product thus reads
  $$
  (\a \d _{\g })(\b \d _{\h }) = \a \abar \g \b \d _{\g \h }.
  \equationmark FormulaProductSimpler
  $$ \endgroup

\section Induction

\label InductionSect As always, we assume the conditions set out in \ref
{StandingHyp}.  From here on the distinguished point $x_0$ mentioned in \ref
{StandingHyp.e} will become important in our development and we will henceforth
use the following notations
  $$\def \quad { }\def \crr {\hfill \cr \pilar {15pt}}
  \matrix {
  \hfill \Sx _{x_0} & := & \{g\in G: x_0\in X_{g\inv }\}, \crr
  \hfill \Hx _{x_0} & := & \{g\in G: x_0\in X_{g\inv },\ \theta _g(x_0)=x_0\},
\crr
  \Orb (x_0) & := & \{\theta _g(x_0): g\in \Sx \}.\crr
  }
  \equationmark BasicNotations
  $$

Whenever there is no question as to which point $x_0$ we are referring to, as it
will often be the case, we will omit subscripts and write $\Sx $ and $\Hx $ in
place of $\Sx _{x_0}$ and $\Hx _{x_0}$, resectively.

Notice that $\Hx $ is a subgroup of $G$, often called the \"{isotropy} group of
$x_0$.  On the other hand observe that $\Sx \Hx \subseteq \Sx $, so $\Sx $ is a
union of left $\Hx $-classes.

The map
  $$
  g\in \Sx \longmapsto \theta _g(x_0) \in \Orb (x_0)
  $$
  is clearly onto, and two elements $g_1$ and $g_2$ in $\Sx $ satisfy $\theta
_{g_1}(x_0)=\theta _{g_2}(x_0)$, if and only if they lie in the same left $\Hx
$-class.

A central ingredient in the induction process to be introduced shortly, is the
subspace $\M $ of the group algebra $\GpAlg G$ given by
  $$
  \M = \hbox {span}\{\delta _g:g\in \Sx \}.
  $$
  As already observed, $\Sx \Hx \subseteq \Sx $, so it follows that $\M $ is
naturally a right $\KH $-module.

Consider the unique bilinear form
  $$
  \absform : \M \times \M \to \KH
  $$
  such that
  $$
  \langle \delta _k,\delta _l\rangle =\left \{
    \matrix {\delta _{k\inv l}, & \hbox {if }k\inv l\in \Hx ,\cr 0, & \hbox
{otherwise}.\hfill \pilar {16pt}}
    \right .
  \equationmark DefineForm
  $$

This may also be written as
  $$
  \langle \delta _k,\delta _l\rangle =\bool {k\inv l\in \Hx } \delta _{k\inv l},
  $$
  where the brackets indicate \"{boolean value}\fn {In fact we shall often use
boolean values in this work, sometimes in a slightly abusive fashion, such as in
  $$
  \bool {x\in X_{g\inv }}f(\theta _g(x)),
  $$
  where $f$ is some scalar valued function on $X$.  The principle behind this is
that, when $x$ is not in the domain $X_{g\inv }$ of $\theta _g$, so that $\theta
_g(x)$ is not defined, the zero boolean value of the expression ``$x\in X_{g\inv
}$" predominates and turns the whole expression into zero.  In other words,
\"{zero} times something which is \"{not defined} is taken to be zero.  It is
true that an excessive abuse of this principle may perhaps lead to unexpected
consequences, but we promisse to use it only to shorten expressions which could
otherwise be writen in two clauses, such as \ref {DefineForm}.}.

  An important property of this form is expressed by the identity
  $$
  \langle m,na\rangle = \langle m,n\rangle a,\for m,n\in \M , \for a\in \KH ,
  \equationmark InnProdRightLinear
  $$
  which the reader may easily prove.

\state Proposition \label QBasis
  If\/ $R\subseteq \Sx $ is a system of representatives of left $\Hx $-classes,
so that
  $$
  \Sx = {\buildrel \thickdot \over {\medcup _{r\in R}}} r\Hx ,
  $$
  then, for all $m$ in $\M $, one has
  $$
  m=\sum _{r\in R}\delta _r\langle \delta _r,m\rangle ,
  $$
  where the sum is always finite in the sense that there are only finitely many
nonzero summands.

  \Proof
  Assuming that $m=\delta _l$, there exists a unique $k$ in $R$ such that $l\Hx
=k\Hx $, which is to say that $k\inv l\in \Hx $.  Then
  $$
  \sum _{r\in R}\delta _r\langle \delta _r,\delta _l\rangle = \sum _{r\in
R}\delta _r\bool {r\inv l\in \Hx }\delta _{r\inv l} = \delta _k\delta _{k\inv l}
= \delta _l.
  $$
  The general case folows by writing $m$ as a linear combination of the $\delta
_l$.
  \endProof

Besides being a right $\KH $-module, $\M $ is also a left module:

\state Proposition
  There is a left ($\alg $)-module structure on $\M $ such that
  $$
  (f\d _g)\delta _l = \bool {gl\in \Sx } f\big (\theta _{gl}(x_0)\big )\delta
_{gl},
  $$
  for every $f\in \lc {X_g}$, and all $l\in \Sx $.  With this, $\M $ moreover
becomes an ($\alg $) - $\KH $ - bimodule.

  \Proof Left for the reader. \endProof

Given any left $\KH $-module $V$, one may therefore build the left $(\alg
)$-module
  $$
  \M \otimes _{\KH }V,
  $$
  henceforth denoted simply by $\M \otimes V$.  This left module structure is
well known, but it might be worth spelling it out here: given $b\in \alg $, one
has
  $$
  b(m\otimes v) = (bm)\otimes v,
  \for m\in \M ,\ v\in V.
  $$

\definition \label DefinEinducedModule
  The left $(\alg )$-module $\M \otimes V$ mentioned above is said to be the
module \"{induced} by $V$.

  Recall that a module $V$ is said to be \"{irreducible} if it has no nontrivial
submodules or, equivalently, if the submodule generated by any nonzero element
coincides with $V$.

\state Proposition \label IndPreserveIrred
  If\/ $V$ is an irreducible\/ $\KH $-module, then $\M \otimes V$ is irreducible
as a $(\alg )$-module.

\Proof
  Given any nonzero vector $w\in \M \otimes V$, we must show that the submodule
it generates, here denoted by $\langle w\rangle $, coincides with $\M \otimes
V$.  In order to do this, write
  $$
  w=\sum _{i=1}^nm_i\otimes u_i,
  $$
  and let $R\subseteq \Sx $ be a system of representatives of left classes for
$\Sx $ modulo $\Hx $.  So by \ref {QBasis} we have
  $$
  w =
  \sum _{i=1}^n\sum _{r\in R}\delta _r\langle \delta _r,m_i\rangle \otimes u_i =
  \sum _{r\in R}\sum _{i=1}^n\delta _r\otimes \langle \delta _r,m_i\rangle u_i =
  \sum _{r\in R}\delta _r\otimes \sum _{i=1}^n\langle \delta _r,m_i\rangle u_i =
  \sum _{r\in R}\delta _r\otimes v_r,
  $$
  where the $v_r$ are defined by the last equality above.  Since all sums
involved are finite, the set
  $$
  \Gamma =\{r\in R:v_r\neq 0\}
  $$
  must be finite.  It is moreover nonempty, since we are assuming that $w\neq
0$.

\def \br {{s}} Fixing any $\br $ in $R$, we claim that $\delta _\br \otimes
v_\br $ lies in $\langle w\rangle $.  To see this, notice that no two elements
of $R$ are in the same left $\Hx $-class, so the points $\theta _r(x_0)$ are
pairwise distinct.  We may then pick some $f$ in $\Lc $ such that $f\big (\theta
_\br (x_0)\big )=1$, while $f\big (\theta _r(x_0)\big )=0$, for all $r\in \Gamma
\setminus \{\br \}$.
  We then have
  $$
  \langle w\rangle \ni fw =
  \sum _{r\in \Gamma }f\delta _r\otimes v_r =
  \sum _{r\in \Gamma }\bool {r\in \Sx }f(\theta _r(x_0))\delta _r\otimes v_r =
  \delta _\br \otimes v_\br .
  $$

  We next show that
  $\langle w\rangle $
  contains $\delta _\br \otimes V$.  Since $V$ is irreducible as a $\KH
$-module, we have that $V$ is spanned by the set
  $
  \{\delta _hv_\br :h\in \Hx \},
  $
  so it is enough to prove that
  $$
  \delta _\br \otimes \delta _hv_\br \in \langle w\rangle , \for h\in \Hx .
  \equationmark DeltaDeltaInMod
  $$

We thus fix some $h$ in $\Hx $, and put $g=shs\inv $.  Observing that
  $$
  \theta _s(x_0) = \theta _s\big (\theta _h(x_0)\big ) = \theta _{sh}(x_0) =
\theta _{gs}(x_0),
  $$
  we see that $x_0\in X_{(gs)\inv }\cap X_{s\inv }$.  Consequently
  $$
  \theta _{gs}(x_0)\in \theta _{gs}\big (X_{(gs)\inv }\cap X_{s\inv }\big ) =
X_{gs}\cap X_g.
  $$
  We may then choose some $f$ in $\lc {X_g}$ such that $f\big (\theta
_{gs}(x_0)\big )=1$, and then
  $$
  \langle w\rangle \ni
  f\d _g(\delta _\br \otimes v_\br ) =
  \bool {gs\in \Sx }f\big (\theta _{gs}(x_0)\big )\delta _{g\br }\otimes v_\br =
  \delta _{g\br }\otimes v_\br =
  \delta _{sh}\otimes v_\br =
  \delta _s\delta _h\otimes v_\br =
  \delta _s\otimes \delta _hv_\br ,
  $$
  thus proving \ref {DeltaDeltaInMod}, and hence that
  $
  \delta _\br \otimes V \subseteq \langle w\rangle .
  $

  We will conclude the proof by showing that $\delta _k\otimes V \subseteq
\langle w\rangle $, for every $k$ in $\Sx $.  Given any such $k$, set $g=ks\inv
$, and notice that
  $
  x_0\in X_{k\inv }\cap X_{s\inv },
  $
  so
  $$
  \theta _k(x_0)\in \theta _k\big (X_{k\inv }\cap X_{\br \inv }\big ) =X_k\cap
X_{k\br \inv } \subseteq X_g.
  $$

  So we may find some $f$ in $\lc {X_g}$ such that $f\big (\theta _k(x_0)\big
)=1$, and for every $v$ in $V$, one has
  $$
  \langle w\rangle \ni
  f\d _g(\delta _\br \otimes v) =
  \bool {g\br \in \Sx }f\big (\theta _{g\br }(x_0)\big )\delta _{g\br }\otimes v
=
  \bool {k\in \Sx }f\big (\theta _{k}(x_0)\big )\delta _{k}\otimes v =
  \delta _{k}\otimes v.
  $$
  This shows that $\delta _k\otimes V \subseteq \langle w\rangle $, as desired,
and hence that $\langle w\rangle = \M \otimes V$, concluding the proof.
\endProof

Our next goal will be to compute the annihilator of the induced module in terms
of the annihilator of the original module $V$.  We begin with a useful technical
result.

\state Lemma \label VectorAnnnihiletor
  Let $V$ be a left $\KH $-module and let $I$ be the annihilator of\/ $V$ in
$\KH $.
  Given $m\in M$, the following are equivalent:
  \izitem
  \zitem $m\otimes v=0$, for all $v$ in $V$,
  \zitem $\langle n,m\rangle \in I$, for all $n\in \M $.

\Proof (ii) $\imply $ (i): Let $R$ be a system of representatives of left $\Hx
$-classes in $\Sx $.  Then for every $v$ in $V$ we have
  $$
  m\otimes v \={QBasis}
  \sum _{r\in R}\delta _r\langle \delta _r,m\rangle \otimes v =
  \sum _{r\in R}\delta _r\otimes \langle \delta _r,m\rangle v = 0.
  $$

\noindent (i) $\imply $ (ii): Fixing $n\in \M $, consider the bilinear mapping
  $$
  (m,v)\in \M \times V \mapsto \langle n,m\rangle v \in V.
  $$
  By \ref {InnProdRightLinear} this is $\KH $-balanced, so there is a well
defined $\Fld $-linear mapping
  $T_n: \M \otimes V \to V$, such that
  $$
  T_n(m\otimes v) = \langle n,m\rangle v,
  $$
  for all $m$ in $\M $, and all $v$ in $V$.
  Assuming that $m$ satisfies (i) we then have that
  $$
  \langle n,m\rangle v=T_n(m\otimes v)=0, \for n\in \M ,\for v\in V,
  $$
  so $\langle n,m\rangle $ annihilates $V$, whence $\langle n,m\rangle \in I$,
proving (ii).
  \endProof

As a consequence we obtain the following description of the annihilator of an
induced module.

\state Corollary \label IndIdealDependsOnlyIdeal
  Let $V$ be a left $\KH $-module and let $I$ be the annihilator of\/ $V$ in
$\KH $.  Then the annihilator of $\M \otimes V$ in $\alg $ is given by
  $$
  \{b\in \alg : \langle n, bm\rangle \in I,\ \forall n,m\in \M \}.
  $$
  In particular the annihilator of $\M \otimes V$ depends only on $I$.

\Proof
  One has that $b$ lies in the annihilator of $\M \otimes V$, iff $bm\otimes
v=0$, for all $m$ and $v$, which is equivalent to saying that
  $\langle n,bm\rangle \in I$, for all $n$ and $m$, by \ref
{VectorAnnnihiletor}.  \endProof

Since the annihilator of $\M \otimes V$ depends only on $I$, rather than on $V$,
we may think of it as built out of $I$.  To account for this we give the
following:

\definition \label IntroJXzeroI
  Given any
  ideal\fn {Unless otherwise stated, all ideals in this paper are assumed to be
two-sided.}
  $I\ideal \KH $, we shall let
  $$
  \JxComXZeroI := \{b\in \alg : \langle n, bm\rangle \in I,\ \forall n,m\in \M
\}.
  $$
  This will be referred to as the \"{ideal induced by $I$}.
  When there is no risk of confusion we shall write this simply as $\JxI $.

Reinterpreting \ref {IndIdealDependsOnlyIdeal} with the terminology just
introduced, we have:

\state Proposition \label AnnInducedVsIndIdeal
  Let $V$ be a left $\KH $-module and let $I$ be the annihilator of\/ $V$ in
$\KH $.  Then the annihilator of $\M \otimes V$ coincides with the ideal induced
by $I$.

So far it is clear that $\JxI $ is a \"{right} ideal, but we will shortly prove
that it is indeed a two-sided ideal.

The behavior of the induction process under inclusion and intersection is easy
to understand:

\state Proposition
  \izitem
  \zitem If $I_1$ and $I_2$ are ideals of\/ $\KH $ with $I_1\subseteq I_2$, then
$\Jx {I_1}\subseteq \Jx {I_2}$.
  \zitem Given any family $\{I_\lambda \}_{\lambda \in \Lambda }$ of ideals of\/
$\KH $, then
  $
  \JxBig {\bigcap _{\lambda \in \Lambda }I_\lambda } = \bigcap _{\lambda \in
\Lambda }\Jx {I_\lambda }.
  $

\Proof Follows easily by inspecting the definitions involved.  \endProof

When checking that $\langle n, bm\rangle \in I$, for all $n,m\in \M $, as
required by the above definition, it suffices to consider $m=\delta _k$ and
$n=\delta _l$, for $k,l\in \Sx $, since these generate $\M $.  It is therefore
nice to have an explicit formula for use in this situation:

\state Proposition \label DeltaKbDeltaL Given $b=\sum _{g\in G}f_g\d _g$ in
$\alg $, and given $k$ and $l$ in $\Sx $, one has that
  $$
  \langle \delta _k, b\delta _l\rangle = \sum _{g\in k\Hx l\inv }f_g\big (\theta
_k(x_0)\big )\, \delta _{k\inv gl}.
  $$

\Proof
  We have
  $$
  \langle \delta _k, b\delta _l\rangle =
  \sum _{g\in G}\langle \delta _k, (f_g\d _g)\delta _l\rangle =
  \sum _{g\in G}\langle \delta _k, \bool {gl\in \Sx }f_g\big (\theta
_{gl}(x_0)\big )\delta _{gl}\rangle \quebra =
  \sum _{g\in k\Hx l\inv }\bool {gl\in \Sx }f_g\big (\theta _{gl}(x_0)\big )\,
\delta _{k\inv gl} = \cdots
  $$

  If $g\in k\Hx l\inv $, then $gl$ lies in $k\Hx $, so $\theta _{gl}(x_0)$ is
indeed defined, meaning that $gl\in \Sx $, and $\theta _{gl}(x_0)$ coincides
with $\theta _k(x_0)$.
  So the above equals
  $$
  \cdots =
  \sum _{g\in k\Hx l\inv }f_g\big (\theta _k(x_0)\big )\, \delta _{k\inv gl},
  $$
  concluding the proof.
  \endProof

The above computation allows for a very concrete criteria for membersip in $\JxI
$, namely:

\state Proposition \label UseDeltaKbDeltaL
  Given any ideal $I\ideal \KH $, and given any $b=\sum _{g\in G}f_g\d _g$ in
$\alg $, one has that $b\in \JxI $, if and only if, for every $k$ and $l$ in
$\Sx $, one has that
  $$
  \sum _{g\in k\Hx l\inv }f_g\big (\theta _k(x_0)\big )\, \delta _{k\inv gl} \in
I.
  $$

\Proof Follows from \ref {DeltaKbDeltaL} and the fact that the $\delta _k$
generate $\M $, as a $\Fld $-vector space.  \endProof

Let us now discuss two trivial examples:

\state Proposition \label Examples
  \iaitem
  \aitem If\/ $I=\KH $, then
  $\JxI $ coincides with the whole algebra $\alg $.
  \aitem If\/ $I=\{0\}$, then
  $$
  \JxI =\Big \{\medsum \nolimits _{g\in G}f_g\d _g\in \alg : f_g|_{\cOrb
(x_0)}=0,\ \forall g\in G\Big \}.
  $$

\Proof
  The first statement is clear.  As for (b), first notice that a locally
constant function vanishing on $\Orb (x_0)$, necessarily also vanishes on the
closure $\cOrb (x_0)$.  This said, let
  $$
  b=\sum _{g\in G}f_g\d _g \in \alg .
  $$
  Assuming that $b$ lies in $\JxI $, and
  given any point $y$ in the orbit of $x_0$, we will prove that $f_g(y)=0$, for
all $g$.
  In case $y\notin X_g$, it is clear that $f_g(y)=0$, since the support of $f_g$
is contained in $X_g$.  Otherwise, if $y\in X_g$, write $y=\theta _k(x_0)$, for
some $k\in \Sx $, and observe that $y\in X_k\cap X_g$, so
  $$
  x_0 = \theta _{k\inv }(y) \in \theta _{k\inv }(X_k\cap X_g) = X_{k\inv }\cap
X_{k\inv g},
  $$
  from where we see that $l:=g\inv k$ lies in $\Sx $.  Consequently
  $$
  \{0\} = I \ni \langle \delta _k,b\delta _l\rangle \={DeltaKbDeltaL}
  \sum _{g'\in k\Hx l\inv }f_{g'}\big (\theta _k(x_0)\big )\, \delta _{k\inv
g'l}.
  $$
  Among the above summands, one is to find
  $
  g'=k1l\inv = kk\inv g = g,
  $
  so in particular
  $$
  0=f_{g'}\big (\theta _k(x_0)\big ) = f_g(y).
  $$
  This shows that $f_g$ vanishes on the orbit of $x_0$, and hence also on its
closure.

Conversely, assuming that each $f_g$ vanishes on the orbit of $x_0$, it is clear
from \ref {DeltaKbDeltaL} that
  $
  \langle \delta _k,b\delta _l\rangle = 0 \in I,
  $
  so $b\in \JxI $.  \endProof

Much has been said about the intersection of an ideal in a crossed product
algebra and its intersection with the coefficient algebra.  In the case of
induced ideals we have:

\state Proposition \label InterCoefAlg
  Let $I$ be a
  proper\fn {We say that an ideal in an algebra is proper when it is not equal
to the whole algebra.}
  ideal in $\KH $. Then the intersection $\JxI \cap \,\Lc $ consists of all $f$
in $\Lc $ vanishing on $\cOrb (x_0)$.

\Proof
  Let $f\in \JxI \cap \Lc $.  Then, choosing any $k$ in $\Sx $, we have
  $$
  I\ni \langle \delta _k,f\delta _k\rangle \={DeltaKbDeltaL}
  f\big (\theta _k(x_0)\big ) \delta _1.
  $$
  Should $f\big (\theta _k(x_0)\big )$ not vanish, the above would be an
invertible element in $I$, whence $I=\KH $, contradicting the hypothesis.  Thus
$f\big (\theta _k(x_0)\big )=0$, showing that $f$ vanishes on the orbit of
$x_0$, and hence also on its closure.

Conversely, if $f$ vanishes on $\cOrb (x_0)$, then by \ref {Examples.ii}
  $$
  f\in \Jx {\{0\}}\subseteq \JxI .
  \closeProof
  $$
  \endProof

Consider the map
  $$
  E=E_{x_0}:\alg \to \cp \Hx ,
  \equationmark IntroE
  $$
  given by
  $$
  E\Big (\sum _{g\in G}f_g\d _g\Big ) = \sum _{h\in \Hx }f_h\d _h.
  $$
  This is sometimes called a \"{conditional expectation}.  One of its important
properties is that of being a $(\cp \Hx )$-bimodule map in the sense that if
$a\in \cp \Hx $, and $b\in \alg $, then
  $$
  E(ab)= aE(b) \and E(ba)= E(b)a.
  \equationmark CondExp
  $$
  It is also evident that $E$ is a projection from $\alg $ onto $\cp \Hx $.

  Consider also the map
  $$
  \nu =\nu _{x_0}:\cp \Hx \to \KH ,
  \equationmark IntroNu
  $$
  given by
  $$
  \nu \Big (\sum _{h\in \Hx }f_h\d _h\Big ) = \sum _{h\in \Hx }f_h(x_0)\delta
_h.
  $$
  Since $x_0$ is fixed by $\Hx $, one may easily show that $\nu $ is an algebra
homomorphism.

The composition of $\nu $ and $E$ is therefore the map
  $$
  F=F_{x_0}=\nu \circ E:\alg \to \KH ,
  \equationmark IntroF
  $$
  given by
  $$
  F\Big (\sum _{g\in G}f_g\d _g\Big ) = \sum _{h\in \Hx }f_h(x_0)\delta _h.
  $$

There is a useful relationship between $F$ and the above bilinear form $\absform
$, expressed as follows:

\state Lemma \label FvsBilForm
  Let $k,l\in G$, and choose $p\in \lc {X_{k\inv }}$, and $q\in \lc {X_l}$, so
that
  defining
  $$
  u=p\d _{k\inv }, \and v=q\d _l,
  $$
  one has that $u$ and $v$ are in $\alg $.
  Then, for every $b$ in $\alg $, one has that
  $$
  F(ubv) =
  \left \{\matrix {
    p(x_0)\,q\big (\theta _l(x_0)\big )\,\langle \delta _k, b\delta _l\rangle ,
& \hbox { if } k,l\in \Sx ,\hfill \stake {15pt} \cr
    0, & \hbox { otherwise}.\hfill
  }\right .
  $$

\Proof
  Write $b=\sum _{g\in G}f_g\d _g$, so that
  $$
  E(ubv) =
  E\Big ( \sum _{g\in G} p\d _{k\inv }\ f_g\d _g\ q\d _l\Big ) =
  E\Big ( \sum _{g\in G} p\,\alpha _{k\inv }(f_g1_k)\,\alpha _{k\inv
g}(q1_{g\inv k})\d _{k\inv gl}\Big ) \quebra =
  \sum _{g\in k\Hx l\inv } p\,\alpha _{k\inv }(f_g1_k)\,\alpha _{k\inv
g}(q1_{g\inv k})\d _{k\inv gl}.
  $$
  Therefore
  $$
  F(ubv) =
  \nu \big (E(ubv)\big ) =
  \sum _{g\in k\Hx l\inv } p(x_0)\ \alpha _{k\inv }(f_g1_k)|_{x_0}\ \alpha
_{k\inv g}(q1_{g\inv k})|_{x_0}\ \delta _{k\inv gl} \quebra =
  \sum _{g\in k\Hx l\inv } p(x_0)\,\bool {x_0\in X_{k\inv }}f_g\big (\theta
_k(x_0)\big )\,\bool {x_0\in X_{k\inv g}}q\big (\theta _{g\inv k}(x_0)\big
)\delta _{k\inv gl} = \cdots
  $$ Notice that, whenever $r\inv s\in \Hx $, one has that $\theta _{r\inv
s}(x_0)=x_0$, so
  $$
  x_0\in X_{r\inv } \iff x_0\in X_{s\inv },
  $$
  and if these equivalent conditions hold then
  $$
  \theta _r(x_0)=\theta _s(x_0).
  $$
  Applying this to $r=g\inv k$, and $s=l$, we see that the above equals
  $$
  \cdots =
  \sum _{g\in k\Hx l\inv } p(x_0)\, \bool {x_0\in X_{k\inv }}f_g\big (\theta
_k(x_0)\big )\, \bool {x_0\in X_{l\inv }}q\big (\theta _l(x_0)\big )\delta
_{k\inv gl} \quebra =
  \bool {x_0\in X_{k\inv }} \bool {x_0\in X_{l\inv }}\, p(x_0)\,q\big (\theta
_l(x_0)\big )\sum _{g\in k\Hx l\inv } f_g\big (\theta _k(x_0)\big )\delta
_{k\inv gl}
  \={DeltaKbDeltaL}
  \bool {k,l\in \Sx }\, p(x_0)\,q\big (\theta _l(x_0)\big )\, \langle \delta _k,
b\delta _l\rangle .
  \closeProof
  $$
  \endProof

Let us now use the above result with the purpose of giving an alternative
definition of $\JxI $, where $F$ is employed instead of the form $\absform $.

\state Proposition \label AltDescription
  Given any ideal $I\ideal \KH $, one has that
  $$
  \JxI = \{b\in \alg : F(ubv)\in I,\ \forall u,v\in \alg \}.
  $$

\Proof
  In order to prove the statement we must show that, for any given $b\in \alg $,
the following are equivalent:
  \izitem
  \zitem $F(ubv)\in I$, for all $u,v\in \alg $,
  \zitem $\langle n, bm\rangle \in I$, for all $n,m\in \M $.

\bigskip \noindent (i) $\imply $ (ii): It is clearly enough to prove (ii) for
$n=\delta _k$, and $m=\delta _l$, where $k$ and $l$ are arbitrary elements of
$\Sx $.  In order to do this, pick $p\in \lc {X_{k\inv }}$, and $q\in \lc
{X_l}$, such that $p(x_0)=1$, and $q\big (\theta _l(x_0)\big )=1$. Then
  $$
  \langle n, bm\rangle =
  \langle \delta _k, b\delta _l\rangle =
  p(x_0)\,q\big (\theta _l(x_0)\big )\, \langle \delta _k, b\delta _l\rangle
\={FvsBilForm}
  F(ubv)\in I,
  $$
  proving (ii).

\bigskip \noindent (ii) $\imply $ (i): It is clearly enough to prove (i) for
$u=p\d _{k\inv }$, and $v=q\d _l$, where $k$ and $l$ are arbitrary elements of
$G$, while $p\in \lc {X_{k\inv }}$, and $q\in \lc {X_l}$.  If $k$ and $l$ lie in
$\Sx $, we have that
  $$
  F(ubv) \={FvsBilForm}
  p(x_0)q\big (\theta _l(x_0)\big )\langle \delta _k, b\delta _l\rangle \in I.
  $$
  On the other hand, if either $k$ or $l$ are not in $\Sx $, then again by \ref
{FvsBilForm}, we have that $F(ubv)=0\in I$.
  \endProof

When $I$ is the annihilator of a left $\KH $-module $V$, we have seen that $\JxI
$ is the annihilator of $\M \otimes V$, hence a two-sided ideal.  Should there
be any doubt that $\JxI $ is always a two-sided ideal (in case $I$ is not
presented\fn {Notice, however, that any ideal of a unital algebra is the
annihilator of some left module, namely the quotient algebra.} as the
annihilator of some left $\KH $-module), the above description of $\JxI $ may be
used to dispell this doubt.

\section Admissible ideals

\label AdmIdSect
  As always we assume \ref {StandingHyp}.  So far we have not considered the
question of which ideals $I\ideal \KH $ actually lead to nontrivial induced
ideals.  In case $\theta $ is topologically free and minimal, a situation well
known to lead to a simple crossed product \cite [Theorem 4.1]{Dan}, at least
some points $x_0$ in $X$ are allowed to possess a nontrivial isotropy group $\Hx
$, and hence there might be plenty ideals $I\ideal \KH $ to choose from, but the
simplicity of $\alg $ prevents induced ideals from being nontrivial.

In this section we shall begin to explore the delicate relationship between
ideals of the isotropy group algebra and the induced ideals they lead to.

The map $F_{x_0}$, introduced in \ref {IntroF}, will play a crucial role in this
section.  Because we will always consider the induction process relative to the
point $x_0$ fixed in \ref {StandingHyp.e}, we will abolish the subscript writing
$F$ in place of $F_{x_0}$.

\state Proposition \label FAlmostMult
  Given $a$ and $b$ in $\alg $, such that either $a$ or $b$ is in $\cp \Hx $,
then
  $$
  F(ab) = F(a)F(b).
  $$

\Proof Suppose that $b$ is in $\cp \Hx $.  Then, using \ref {CondExp}, we have
  $$
  F(ab) =
  \nu \big (E(ab)\big ) =
  \nu \big (E(a)b\big ) =
  \nu \big (E(a)E(b)\big ) =
  \nu \big (E(a)\big ) \nu \big (E(b)\big ) =
  F(a) F(b).
  $$
  A similar reasoning applies when $a$ is in $\cp \Hx $.  \endProof

In particular, if $\varphi $ is in $\Lc $, then
  $$
  F(a\varphi ) = \varphi (x_0)F(a) = F(\varphi a), \for a\in \alg .
  \equationmark FunctOut
  $$
  Another instance of \ref {FAlmostMult} is obtained when $\varphi $ is
supported on $X_h$, for some $h$ in $\Hx $, in which case we have
  $$
  F(a\,\varphi \d _h) = F(a)\varphi (x_0)\delta _h, \and F(\varphi \d _h\,a) =
\varphi (x_0)\delta _h\,F(a).
  \equationmark FunctOutTwo
  $$

\state Proposition \label MapIdeals
  If $J$ is any ideal in $\alg $, then $F(J)$ is an ideal in $\KH $.

\Proof
  Given $c\in F(J)$, and $d\in \KH $, we must prove that $cd$ and $dc$ lie in
$F(J)$, and it clearly suffices to assume that $d=\delta _h$, for some $h$ in
$\Hx $.  Choose $b$ in $J$ such that $F(b)=c$, and let $\varphi \in \lc {X_h}$
be such that $\varphi (x_0)=1$.  Then
  $$
  cd =
  c\varphi (x_0)\delta _h =
  F(b)F(\varphi \d _h) \= {FunctOutTwo}
  F(b\varphi \d _h) \in F(J).
  $$
  A similar reasoning proves that $dc\in F(J)$.
  \endProof

Applying this to induced ideals we get the following:

\state Proposition \label IntroIPrime
  Let $I$ be an ideal in $\KH $, and put
  $
  I'= F(\JxI ).
  $
  Then
  \izitem
  \zitem $I'$ is an ideal of $\KH $,
  \zitem $I'\subseteq I$,
  \zitem $\Jx {I'} = \JxI $.

\Proof
  (i) Follows from \ref {MapIdeals}.

\medskip \noindent (ii) Given $a$ in $I'$, write $a=F(b)$, with $b\in \JxI $.
Choosing $u\in \Lc $, with $u(x_0)=1$, we then have
  $$
  I \explain {AltDescription}\ni
  F(ubu) \={FunctOut}
  u(x_0)F(b)u(x_0) = F(b) = a.
  $$
  This proves (ii).

\medskip \noindent (iii) Since $I'\subseteq I$, it is obvious that $\Jx
{I'}\subseteq \JxI $.  On the other hand, if $b\in \JxI $, then for every $u$
and $v$ in $\alg $, one has that $ubv\in \JxI $, whence
  $$
  F(ubv)\in F(\JxI ) = I'.
  $$
  This proves that $b\in \Jx {I'}$, concluding the proof.  \endProof

The grand conclusion of this result is that, should we want to catalogue all
induced ideals, we do not need to consider all ideals $I\ideal \KH $, since $I$
may be replaced by $I'$, without affecting the outcome of the induction process.

This motivates the question of how to separate the ideals that matter from those
which don't, a task we now begin to undertake.

\definition \label DefAdmiss An ideal $I\ideal \KH $ is said to be
\"{admissible} if $F(\JxI )=I$.

Interpreting \ref {IntroIPrime} from the point of view of the concept just
introduced we have:

\state Corollary \label ChangeIdeals
  For every ideal $I\ideal \KH $, there exists a unique admissible ideal
$I'\subseteq I$, such that $\JxI =\Jx {I'}$.

\Proof
  Set $I'=F(\JxI )$.  Then $\JxI =\Jx {I'}$, by \ref {IntroIPrime.iii}.
Moreover
  $$
  F(\Jx {I'})=F(\JxI ) = I',
  $$
  so $I'$ is admissible.
  If $I'$ and $I''$ are two admissible ideals inducing the same ideal of $\alg
$, then
  $$
  I'=F(\Jx {I'})=F(\Jx {I''})=I''.
  \closeProof
  $$
  \endProof

We have in fact already encountered examples of admissible ideals:

\state Proposition \label TrivialIdeals
  The two trivial ideals of\/ $\KH $, namely $\{0\}$ and $\KH $, itself, are
admissible.

\Proof
  Setting $I=\{0\}$, we have that
  $$
  I=\{0\}\subseteq F(\JxI ) \explain {IntroIPrime.ii}\subseteq I,
  $$
  so $I$ is admissible.

On the other hand, if $I=\KH $, we have seen in \ref {Examples} that $\JxI =\alg
$, so
  $$
  F(\JxI )=F\big (\alg \big )=\KH =I,
  $$
  so, again, $I$ is seen to be admissible.
  \endProof

In order to better understand admissible ideals, we must be able to describe the
image of an ideal in $\alg $ through $F$.

\state Proposition \label LiftC
  Let $J\ideal \alg $ be any ideal and let
  $
  c = \sum _{h\in \Hx }c_h\delta _h
  $
  be any element of\/ $\KH $.  Then $c$ is in $F(J)$ if and only if there exists
a {\klopen } set $V$, such that
  $$
  x_0\in V\subseteq X_h,
  $$
  whenever $c_h\neq 0$, satisfying
  $$
  c_V:= \sum _{h\in \Hx }c_h1_V\d _h \in J.
  $$

\Proof
  Assuming that $c_V$ is in $J$, we have
  $$
  F(J) \ni F(c_V) = \sum _{h\in \Hx }c_h1_V(x_0)\delta _h = \sum _{h\in \Hx
}c_h\delta _h = c,
  $$
  so $c\in F(J)$.  Conversely, if $c\in F(J)$, pick $b$ in $J$ such that
$c=F(b)$.  We will initially prove that $b$ may be chosen in $\cp \Hx $.

Write $b=\sum _{g\in \Gamma }f_g\d _g$, where $\Gamma $ is a finite subset of
$G$, and set
  $$
  \def \quad {\ }\def \crr {\hfill \stake {12pt}\cr }
  \matrix {
  \Gamma _1 &=& \{g\in \Gamma : x_0\notin X_{g\inv }\}, \crr
  \Gamma _2 &=& \{g\in \Gamma : x_0\in X_{g\inv },\ \theta _g(x_0)\neq x_0\},
\crr
  \Gamma _3 &=& \{g\in \Gamma : x_0\in X_{g\inv },\ \theta _g(x_0)=x_0\} =
\Gamma \cap \Hx . \crr  }
  $$
  It is then clear that $\Gamma $ is the disjoint union of the $\Gamma _i$.

For each $g$ in $\Gamma _1$, using that $X_{g\inv }$ is closed, choose a
{\klopen } set $W_g$, such that
  $$
  x_0\in W_g\subseteq X\setminus X_{g\inv }.
  $$

For each $g$ in $\Gamma _2$, choose open sets $U$ and $V$, such that
  $x_0\in U$, $\theta _g(x_0)\in V$, and $U\cap V=\emptyset $.
  By replacing $V$ with $V\cap X_g$ we may assume that $V\subseteq X_g$.  We
then set $Z=U\cap \,\theta _{g\inv }(V)$, and observe that
  $x_0\in Z\subseteq X_{g\inv }$, and that
  $$
  Z\cap \theta _g(Z) \subseteq U\cap \theta _g\big (\theta _{g\inv }(V)\big ) =
U\cap V = \emptyset .
  $$
  Choosing a {\klopen } neighborhood $W_g$ of $x_0$ contained in $Z$, we then
have that
  $$
  x_0\in W_g\subseteq X_{g\inv }, \and W_g\cap \theta _g(W_g) = \emptyset .
  $$

Ignoring $\Gamma _3$ for the time being we
  put
  $$
  W = \bigcap _{g\in \Gamma _1\cup \Gamma _2} W_g,
  $$
  and observe that
  $$
  \def \crr {\hfill \stake {12pt}\cr }
  \matrix {
    W\cap X_{g\inv }=\emptyset , & \for g\in \Gamma _1, & \hbox { and}\crr
    x_0\in W\subseteq X_{g\inv }, \hbox { \ and \ } W\cap \theta _g(W) =
\emptyset , & \for g\in \Gamma _2.}
  $$
  We then have that
  $$
  1_Wb1_W =
  \sum _{g\in \Gamma }1_W(f_g\d _g)1_W =
  \sum _{g\in \Gamma }1_Wf_g\alpha _g(1_W1_{X_{g\inv }})\d _g =
  \sum _{g\in \Gamma }1_Wf_g\alpha _g(1_{W\cap X_{g\inv }})\d _g \quebra =
  \sum _{g\in \Gamma }1_Wf_g1_{\theta _g(W\cap X_{g\inv })}\d _g =
  \sum _{g\in \Gamma }f_g1_{W\cap \theta _g(W\cap X_{g\inv })}\d _g.
  $$

  For $g\in \Gamma _1\cup \Gamma _2$ we have that $W\cap \theta _g(W\cap
X_{g\inv })=\emptyset $, so the summand corresponding to $g$ in the above sum
vanishes.  Therefore,
  $$
  1_Wb1_W =
  \sum _{g\in \Gamma _3}f_g1_{W\cap \theta _g(W\cap X_{g\inv })}\d _g =
  \sum _{g\in \Gamma _3}f'_g\d _g,
  $$
  where the $f'_g$ are defined by the last equality above.  Setting
$b'=1_Wb1_W$, we then have that $b'\in \cp \Hx $, because $\Gamma _3=\Gamma \cap
\Hx $, and moreover $b'\in J$. Recalling that $c=F(b)$, we also have that
  $$
  F(b') =
  F(1_Wb1_W) \={FunctOut}
  1_W(x_0)F(b)1_W(x_0) =
  F(b) = c.
  $$

  Replacing $b$ by $b'$ we have therefore proven our claim that $b$ may be
chosen in $\cp \Hx $, so we are allowed to write
  $$
  b = \sum _{h\in \Hx }f_h\d _h.
  $$

  For each $h$ in $\Hx $, choose a {\klopen } set $V_h\subseteq X_h$ such that
$x_0\in V_h$, and such that $f_h$ is constant on $V_h$.  Letting $V$ be the
intersection of the finitely many $V_h$ for which $f_h$ is nonzero, we have that
the $f_h$ are constant on $V$, so that $1_Vf_h = d_h1_V$, where $d_h$ is the
constant value attained by $f_h$ on V.

We may then define
  $$
  b'':=1_Vb = \sum _{h\in \Hx }1_Vf_h\d _h = \sum _{h\in \Hx }d_h1_V\d _h,
  $$
  observing that, as above, $b''\in J$, and $F(b'')=c$.  The latter may be
expressed as
  $$
  \sum _{h\in \Hx }c_h\delta _h =
  F\big (\sum _{h\in \Hx }d_h1_V\d _h\big ) =
  \sum _{h\in \Hx }d_h1_V(x_0)\delta _h =
  \sum _{h\in \Hx }d_h\delta _h.
  $$
  It follows that $d_h=c_h$, for all $h$, whence
  $$
  c_V = \sum _{h\in \Hx }c_h1_V\d _h = \sum _{h\in \Hx }d_h1_V\d _h = b'' \in J.
  \closeProof
  $$
  \endProof

We may now employ \ref {LiftC} to give a characterization of admissible ideals.

\state Proposition \label AdmissibleCharact
  An ideal $I\ideal \KH $ is admissible if and only if, for every $c=\sum _{h\in
\Hx }c_h\delta _h$ in $I$, there exists a neighborhood\/ $V$ of $x_0$, such that
  $$
  \delta _{k\inv }\Big (\sum _{h\in \Hx \cap k\Hx l\inv }c_h\delta _h\Big
)\delta _l \in I,
  $$
  for all $k$ and $l$ in $\Sx $, such that $\theta _k(x_0)\in V$.

\Proof
  Supposing that $I$ is admissible, pick $c=\sum _{h\in \Hx }c_h\delta _h\in I$.
By hypothesis $c$ is in $F\big (\JxI \big )$, so \ref {LiftC} provides a
{\klopen } set $V\ni x_0$, such that
  $$
  c_V:= \sum _{h\in \Hx }c_h1_V\d _h \in \JxI .
  $$

  In view of the definition of $\JxI $, one has that $\langle \delta
_k,c_V\delta _l\rangle \in I$, for every $k$ and $l$ in $\Sx $, and if we use
\ref {DeltaKbDeltaL} under the hypothesis that $\theta _k(x_0)\in V$, we deduce
that
  $$
  I\ni \langle \delta _k,c_V\delta _l\rangle = \sum _{h\in \Hx \cap k\Hx l\inv
}c_h1_V\big (\theta _k(x_0)\big )\, \delta _{k\inv hl} =
  \sum _{h\in \Hx \cap k\Hx l\inv }c_h\, \delta _{k\inv hl},
  $$
  proving the condition displayed in the statement.

Conversely, assuming that this condition holds, let us show that $I$ is
admissible, namely that $I\subseteq F(\JxI )$, since the reverse inclusion is
granted by \ref {IntroIPrime.ii}.  For this, given $c=\sum _{h\in \Hx }c_h\delta
_h\in I$, pick $V$ as in the statement.  By shrinking $V$ a bit, if necessary,
we may assume that $V$ is {\klopen } and $V\subseteq X_h$, whenever $c_h\neq 0$,
so that
  $$
  c_V:= \sum _{h\in \Hx }c_h1_V\d _h
  $$
  is a legitimate element of $\alg $.  We then claim that $c_V\in \JxI $.  To
prove this it is enough to verify that
  $\langle \delta _k,c_V\delta _l\rangle \in I$, for all $k$ and $l$ in $\Sx $.
  Using \ref {DeltaKbDeltaL} again we have
  $$
  \langle \delta _k,c_V\delta _l\rangle =
  \sum _{h\in \Hx \cap k\Hx l\inv }c_h1_V\big (\theta _k(x_0)\big )\, \delta
_{k\inv hl} =
  \bool {\theta _k(x_0)\in V}\, \delta _{k\inv }\Big (\sum _{h\in \Hx \cap k\Hx
l\inv }c_h\, \delta _h\Big )\delta _l.
  $$
  In case $\theta _k(x_0)\in V$, the hypothesis implies that the above belongs
to $I$, and otherwise $\langle \delta _k,c_V\delta _l\rangle $ vanishes so it
also lies in $I$.  This shows that $c_V$ is in $\JxI $, so
  $$
  c = F(c_V) \in F(\JxI ),
  $$
  concluding the proof.  \endProof

Recall from \ref {IntroIPrime.ii} that, for every ideal $I\ideal \KH $, one has
that $F(\JxI )\subseteq I$.  One may similarly inquire about the relationship
between $J$ and $\Jx {F(J)}$.  The answer is given in our next:

\state Proposition \label LargestIdeal
  \izitem
  \zitem For every ideal $J\ideal \alg $, one has that $J\subseteq \Jx {F(J)}$.
  \zitem For every ideal $I\ideal \KH $, one has that $\JxI $ is the largest
among the ideals $J\ideal \alg $ satisfying $F(J)\subseteq I$.

\Proof
  (i) Given $b$ in $J$, notice that for every $u$ and $v$ in $\alg $, one has
that $ubv\in J$, so
  $$
  F(ubv)\in F(J).
  $$
  We then deduce from \ref {AltDescription} that $b$ lies in $\Jx {F(J)}$. This
proves (i).

\medskip \noindent (ii)
  As already mentioned, \ref {IntroIPrime.ii} gives $F(\JxI )\subseteq I$, so
$\JxI $ is indeed among the ideals mentioned above.  Next,
  given any ideal
  $
  J\ideal \alg ,
  $
  with $F(J)\subseteq I$, we have
  $$
  J\explica {(i)}\subseteq \Jx {F(J)}\subseteq \JxI .
  \closeProof
  $$ \endProof

Our main interest is to construct ideals in $\alg $ from admissible ideals in
$\KH $, but it is interesting to remark that one may also go the other way:

\state Proposition \label FJIsAdmissible
  Let $J\ideal \alg $ be any ideal.  Then $F(J)$ is an admissible ideal of\/
$\KH $.

\Proof
  By \ref {LargestIdeal.i} we have that $J\subseteq \Jx {F(J)}$.   So, if we
apply $F$ on both sides of this inclusion, we get
  $$
  F(J)\subseteq F\big (\Jx {F(J)}\big )\explain {IntroIPrime.ii}\subseteq F(J),
  $$
  so we see that $F\big (\Jx {F(J)}\big )=F(J)$, which is to say that $F(J)$ is
admissible.  \endProof

\def \linsp #1{\left [#1\right ]}  \def \w {\eta } \def \prodx
{\prod _{x\in X}} \def \sumx {\bigoplus _{x\in X}} \def \th {{\hbox {\eightrm
th}}} \def \gl {\hbox {GL}} \def \cCHb {b}

\section Representations

\label RepSec As before we adopt our standing assumptions \ref {StandingHyp}.
In this section we will begin the preparations for proving that any
  ideal (always meaning two-sided ideal)
  of $\alg $ is the intersection of ideals induced from isotropy subgroups.

Our methods will largely rely on representation theory, so we begin by spelling
out a trivial connection between representations and ideals.

\state Proposition \label IdealIsKernel
  Let $B$ be a $\Fld $-algebra possessing \"{local units}\fn {Recall that $B$ is
said to have local units if, for every $b$ in $B$, there exists an idempotent
$e\in B$, such that $eb=b=be$.}
  and let us be given an ideal $J\ideal B$.  Then there exists a vector space
$V$, and a
  non-degenerate\fn {We say that $\pi $ is non-degenerate if $V=\linsp {\pi
(B)V}$, brackets meaning linear span.} representation
  $$
  \pi :B\to L(V),
  $$
  such that $J=\Ker (\pi )$.

\Proof
  Let $V=B/J$, denote the quotient map by $q:B\to V$, and consider the
representation $\pi :B\to L(V)$ given by
  $$
  \pi (b)q(v) = q(bv), \for b,v\in B.
  $$
  It is then obvious that $J\subseteq \Ker (\pi )$, but the reverse inclusion
may also be verified: in fact, if $b$ is in $\Ker (\pi )$, choose $e$ in $B$
such that $b=be$, so
  $$
  0=\pi (b)q(e) = q(be) = q(b),
  $$
  whence $b$ is in $J$.  In order to show that $\pi $ is non-degenerate, pick
any $\xi $ in $V$, and write $\xi =q(b)$, for some $b$ in $B$.  Letting $e$ be
such that $b=eb$, we have
  $$
  \xi =q(b)=q(eb)=\pi (e)q(b)\in \pi (B)V.
  \closeProof
  $$ \endProof

To see that the above result applies to our situation, we give the following:

\state Proposition \label LocalUnits
  For every $b$ in $\alg $, there is an idempotent $e\in \Lc $, such that
$eb=b=be$.  In particular $\alg $ has local units.

\Proof
  Given $b$ in $\alg $, write
  $$
  b=\sum _{g\in G\Gamma }f_g\d _g,
  $$
  where $\Gamma $ is a finite subset of $G$, and each $f_g$ lies in $\lc {X_g}$.
For each $g$ in $\Gamma $, let $C_g=\supp (f_g)$, so that $C_g$ is a {\klopen }
subset of $X$, contained in $X_g$.
  Put
  $$
  C=\bigcup _{g\in \Gamma } \Big (C_g\cup \theta _{g\inv }(C_g)\Big ).
  $$

  It follows that $C$ is also a {\klopen } subset of $X$, whence its
characteristic function $1_C$ is an idempotent element of $\Lc $, hence also of
$\alg $.

We next claim that $1_Cb=b1_C=b$.  To see this we observe that, for obvious
reasons, $1_Cf_g=f_g$, for any $g$ in $\Gamma $, so it is clear that $1_Cb=b$.
On the other hand
  $$
  b1_C=
  \sum _{g\in G\Gamma }f_g\d _g1_C \={FormulaProductSimpler} \sum _{g\in G\Gamma
}f_g\abar g{1_C}\d _g,
  \equationmark BOneC
  $$
  while, for every $g$ in $\Gamma $, we have that
  $$
  f_g\abar g{1_C} =
  1_{C_{g}}f_g\alpha (1_C1_{g\inv }) =
  1_{C_{g}}f_g1_{\theta _g({C\cap X_{g\inv }}))} =
  f_g1_{C_g\cap \theta _g({C\cap X_{g\inv }}))} = f_g1_{C_g} = f_g,
  $$
  where, in the penultimate step we have used that
  $$
  C_g =
  \theta _g\big (\theta _{g\inv }(C_g)\big ) =
  \theta _g\big (\theta _{g\inv }(C_g)\cap X_{g\inv }\big ) \subseteq
  \theta _g\big (C\cap X_{g\inv }\big ).
  $$

  This proves that $f_g\abar g{1_C}=f_g$, whence the computation in \ref {BOneC}
gives that $b1_C=b$.  \endProof

From this point on, we will fix an arbitrary ideal $J\ideal \alg $, which in
view of \ref {IdealIsKernel} and \ref {LocalUnits}, we may assume is the kernel
of a likewise fixed non-degenerate representation
  $$
  \pi :\alg \to \lin (V).
  $$

For the time being we will forget about the ideal $J$ mentioned above, and we
will mostly focus our attention on the representation $\pi $, even though our
main long term goal is to study $J$.

\state Proposition \label LcNonDeg
  Regarding the above representation $\pi $, its restriction to $\Lc $ is
non-degenerate.

\Proof Given any vector $\xi $ in $V$, write
  $$
  \xi =\sum _{i=1}^n\pi (b_i)\xi _i ,
  $$
  with $b_i$ in $\alg $, and $\xi _i$ in $V$.   Using \ref {LocalUnits}, for
each $i$ we choose an idempotent $e_i\in \Lc $ such that $b_i=e_ib_i$, so
  $$
  \xi =
  \sum _{i=1}^n\pi (b_i)\xi _i =
  \sum _{i=1}^n\pi (e_ib_i)\xi _i =
  \sum _{i=1}^n\pi (e_i)\pi (b_i)\xi _i \in
  \linsp {\pi \big (\Lc \big )V}.
  \closeProof
  $$ \endProof

\state Proposition \label Disintegration (Disintegration)
  There exists a partial representation
  $$
  u:G\to \lin (V),
  $$
  such that, for all $g\in G$, and $f\in \lc {X}$, one has
  \izitem
  \zitem $u_g\pi (f) = \pi \big (\abar gf\big ) u_g$, \zitemmark CovarCond
  \zitem $\pi (f\Delta _g) = \pi (f)u_g$, provided $f\in \lc {X_g}$, \zitemmark
Reconstruct
  \zitem $u_gu_{g\inv }\pi (f) = \pi (f1_g) = \pi (f)u_gu_{g\inv }$. \zitemmark
Projecao

\Proof Given any $\xi $ in $V$, use \ref {LcNonDeg} to write
  $$
  \xi = \sum _{i=1}^n \pi (\varphi _i) \eta _i,
  $$
  where each $\varphi _i\in \Lc $, and $\eta _i\in V$.  We then define
  $$
  u_g\xi =
  \sum _{i=1}^n \pi \big (\abar g{\varphi _i}\Delta _g \big ) \eta _i.
  \equationmark DefineUg
  $$

  To prove that this is well defined, suppose that $\xi =0$, and let
  $$
  V=\bigcup _{i=1}^n\supp (\varphi _i)\cap X_{g\inv }.
  $$
  So $V$ is a {\klopen } set and $1_V\varphi _i1_{g\inv } = \varphi _i1_{g\inv
}$, for all $i=1,\ldots ,n$.
  Therefore
  $$
  1_{\theta _g(V)}\Delta _g \, \varphi _i =
  \alpha _g(1_V)\abar g{\varphi _i}\Delta _g =
  \abar g{1_V\varphi _i}\Delta _g =
  \abar g{\varphi _i}\Delta _g ,
  $$
  so the right-hand-side of \ref {DefineUg} coincides with
  $$
  \sum _{i=1}^n \pi \big ( 1_{\theta _g(V)}\Delta _g \,\varphi _i\big ) \eta _i
=
  \pi \big (1_{\theta _g(V)}\Delta _g\big ) \sum _{i=1}^n \pi (\varphi _i\big )
\eta _i =
  \pi \big (1_{\theta _g(V)}\Delta _g\big ) \xi = 0.
  $$
  This shows that $u_g$ is well defined.

  In order to prove \ref {LocalCovarCond}, consider a vector $\xi \in V$ of the
form $\xi =\pi (\varphi )\eta $, for some $\varphi \in \Lc $, and $\xi \eta \in
V$, and observe that
  $$
  u_g\pi (f)\xi =
  u_g\pi (f)\pi (\varphi )\eta =
  u_g\pi (f\varphi )\eta =
  \pi \big (\abar g{f\varphi }\Delta _g \big )\eta \quebra =
  \pi \big (\abar gf\big )\pi \big (\abar g\varphi \Delta _g \big )\eta =
  \pi \big (\abar gf\big )u_g\xi .
  $$
  Since the set of vectors $\xi $ of the above form spans $V$, we have proved
\ref {LocalCovarCond}.

  With the goal of proving \ref {LocalReconstruct}, let $\xi =\pi (\varphi )\eta
$, as above, and notice that
  $$
  \pi (f\Delta _g)\xi =
  \pi (f\Delta _g)\pi (\varphi )\eta =
  \pi (f\Delta _g\,\varphi )\eta \quebra =
  \pi \big (f\abar g\varphi \Delta _g \big )\eta =
  \pi (f)\pi \big (\abar g\varphi \Delta _g \big )\eta =
  \pi (f)u_g \xi .
  $$

In order to prove \ref {LocalProjecao} write $f=f_1f_2$, with $f_1,f_2\in \Lc $,
and let $\xi \in V$.  Then
  $$
  u_gu_{g\inv }\pi (f)\xi =
  u_g\pi \big (\abar {g\inv }{f}\d _{g\inv }\big )\xi =
  u_g\pi \big (\abar {g\inv }{f_1}\big )\pi \big (\abar {g\inv }{f_2}\d _{g\inv
}\big )\xi \quebra =
  \pi \big (\abar g{\abar {g\inv }{f_1}}\d _g\big )\pi \big (\abar {g\inv
}{f_2}\d _{g\inv }\big )\xi =
  \pi \big (f_11_g\d _g\,\abar {g\inv }{f_2}\d _{g\inv }\big )\xi \quebra =
  \pi \big (f_11_g\abar g{\abar {g\inv }{f_2}}\big )\xi =
  \pi \big (f_11_gf_21_g\big )\xi =
  \pi \big (f1_g\big )\xi .
  $$

  This proves the first identity in \ref {LocalProjecao}.  As for the second,
let $\xi =\pi (\varphi )\eta $, so
  $$
  \pi (f)u_gu_{g\inv }\xi =
  \pi (f)u_gu_{g\inv }\pi (\varphi )\eta =
  \pi (f)\pi (\varphi 1_g)\eta =
  \pi (f\varphi 1_g)\eta =
  \pi (f1_g)\pi (\varphi )\eta =
  \pi (f1_g)\xi .
  $$
  This concludes the proof of \ref {LocalProjecao}.

We leave it for the reader to prove that $u$ is a partial representation.
\endProof

We are now about to take a major step in our quest to understand general ideals
in terms of induced ones.  Observe that the very definition of induced ideals
requires that a point of $X$ be chosen in advance, so we must begin to see our
representation $\pi $ from the point of view of a chosen point in $X$, a process
which will eventually lead to a \"{discretization} of $\pi $ (see \ref
{Disintegration} below).  This will be accomplished by means of the following
device: for each $x$ in $X$, let
  $$
  I_x = \{f\in \Lc : f(x)=0\},
  $$
  which is clearly an ideal in $\Lc $.  Consequently
  $$
  Z_x:= \linsp {\pi (I_x)V}
  $$
  is invariant under $\Lc $, so there is a well defined representation $\pi _x$
of $\Lc $ on
  $$
  V_x := V/Z_x,
  $$
  making the following diagram commute
  $$\def \quad {}
  \matrix {
  \stake {10pt}
  V & \Flecha {45pt}{\pi (f)} & V \cr
  q_x \Big \downarrow \kern 10pt && \kern 10pt\Big \downarrow q_x \cr
  V_x & \Flecha {45pt}{\pi _x(f)} & V_x
  }
  $$

The first indication that our localization process is bearing fruit is as
follows:

\state Proposition \label EigenVec
  Given any $x$ in $X$, and any $f$ in $\Lc $, one has that
  $$
  \pi _x(f)\eta =f(x)\eta ,
  $$
  for every $\eta $ in $V_x$.

\Proof Using \ref {LcNonDeg} it is enough to verify the statement for $\eta $ of
the form $\eta = q_x\big (\pi (\varphi )\xi \big )$, where $\varphi \in \Lc $,
and $\xi \in V$.  Choose a {\klopen } set $C$ containing $\supp (\varphi )\cup
\{x\}$, and observe that $1_C\varphi =\varphi $.  In addition $f-f(x)1_C$ lies
in $I_x$, since it clearly vanishes at $x$.  Therefore
  $$
  \pi _x(f)\eta =
  \pi _x(f)q_x\big (\pi (\varphi )\xi \big ) =
  q_x\big (\pi (f)\pi (\varphi )\xi \big ) =
  q_x\big (\pi (f\varphi )\xi \big ) = \cdots
  \equationmark qxpifphiv
  $$

  Notice that
  $$
  f\varphi =
  \big (f(x)1_C + f-f(x)1_C\big )\varphi =
  f(x)1_C\varphi + \big (f-f(x)1_C\big )\varphi \explica {(mod $\scriptstyle
I_x$)}\equiv
  f(x)\varphi ,
  $$
  so
  $$
  \pi (f\varphi )\xi \explica {(mod $\scriptstyle Z_x$)}\equiv \pi \big
(f(x)\varphi \big )\xi = f(x)\pi (\varphi )\xi ,
  $$
  and we then conclude that \ref {qxpifphiv} equals
  $$
  \cdots = q_x\big ( f(x)\pi (\varphi )\xi \big ) = f(x) q_x\big ( \pi (\varphi
)\xi \big ) = f(x) \eta .
  \closeProof
  $$
  \endProof

Putting together the definition of $\pi _x$ with the result above, we get the
following useful formulas:
  $$
  q_x\big (\pi (f)\xi \big ) = \pi _x(f) q_x(\xi ) = f(x) q_x(\xi ),
  \equationmark SummarizePix
  $$
  for all $x\in X$, $f\in \Lc $, and $\xi \in V$.

Having focused on $\Lc $, we momentarily lost track of the $u_g$, but there is
still time to bring them back into focus:

\state Proposition
  If $x$ is in $X_{g\inv }$, then:
  \izitem
  \zitem $u_g(Z_x)\subseteq Z_{\theta _g(x)}$, where $u$ is as in \ref
{Disintegration},
  \zitem there exists a linear mapping
  $$
  \mu ^x_g:V_x\to V_{\theta _g(x)},
  $$
  such that
  $$
  \mu ^x_g\big (q_x(\xi )\big ) = q_{\theta _g(x)}(u_g\xi ), \for \xi \in V.
  $$

\Proof (i)
  Let $\xi $ be a vector in $Z_x$ of the form $\xi =\pi (\varphi )\eta $, where
$\varphi \in I_x$, and $\eta \in V$.  Then
  $$
  u_g\xi =
  u_g\pi (\varphi )\eta \={CovarCond}
  \pi \big (\abar g\varphi \big )u_g\eta .
  $$

  Notice that $\abar g\varphi $ lies in $I_{\theta _g(x)}$, because
  $$
  \abar g\varphi |_{\theta _g(x)} =
  \varphi \big (\theta _{g\inv }(\theta _g(x))\big ) =
  \varphi (x)=
  0,
  $$
  whence $u_g\xi \in Z_{\theta _g(x)}$.

\medskip \noindent (ii) Follows immediately from (i).  \endProof

The $\mu _g^x$ obey the following functorial property:

\state Proposition \label MuOnEg
  If $x\in X_{g\inv }\cap X_{g\inv h\inv }$, then the composition
  $$
  V_x \flecha {\mu ^x_g} V_{\theta _g(x)} \Flecha {40pt}{\mu ^{\theta _g(x)}_h}
V_{\theta _{hg}(x)}
  $$
  coincides with $\mu _{hg}^x$.

\Proof We initially claim that for all $g$ in $G$, if $x\in X_g$, then
  $$
  q_x(\xi ) =q_x(u_gu_{g\inv }\xi ), \for \xi \in V.
  $$

  This will clearly follow should we prove that
  $$
  \xi -u_gu_{g\inv }\xi \in Z_x, \for \xi \in V,
  $$
  which we will now do.  By \ref {LcNonDeg}, we may assume that $\xi = \pi
(\varphi )\eta $, for some $\varphi $ in $\Lc $, and $\eta $ in $V$.  We then
have
  $$
  \xi -u_gu_{g\inv }\xi =
  \pi (\varphi )\eta -u_gu_{g\inv }\pi (\varphi )\eta \= {Projecao}
  \pi (\varphi )\eta -\pi (\varphi 1_g)\eta =
  \pi (\varphi -\varphi 1_g)\eta .
  $$
  Observing that $\varphi -\varphi 1_g$ is in $I_x$, we have that $\pi (\varphi
-\varphi 1_g)\eta $ lies in $Z_x$, proving the claim.

Addressing the statement, choose any element of $V_x$, say $q_x(\xi )$, for some
$\xi $ in $V$, and notice that
  $$
  \mu ^{\theta _g(x)}_h\big (\mu _g^x\big (q_x(\xi )\big ) =
  \mu ^{\theta _g(x)}_h\big (q_{\theta _g(x)}(u_g\xi )\big ) =
  q_{\theta _h(\theta _g(x))}(u_hu_g\xi ) \quebra =
  q_{\theta _{hg}(x)}(u_hu_{h\inv }u_hu_g\xi ) =
  q_{\theta _{hg}(x)}(u_hu_{h\inv }u_{hg}\xi ) = \cdots
  $$

Since $\theta _{hg}(x)=\theta _h(\theta _g(x))\in X_h$, and thanks to our claim,
the above equals
  $$
  \cdots \= {MuOnEg}
  q_{\theta _{hg}(x)}(u_{hg}\xi )=
  \mu _{hg}^x(q_x\xi ).
  \closeProof
  $$ \endProof

Let us now consider the representation
  $$
  \Pi = \prodx \pi _x
  $$ of $\Lc $ on the cartesian product $\prodx V_x$.  Thus, if $f\in \Lc $, and
$\w = (\w _x)_{x\in X}\in \prodx V_x$, we have
  $$
  \big (\Pi (f)\w \big )_x = \pi _x(f)\w _x,
  \for x\in X.
  $$

  Incidentally, by \ref {EigenVec} the term $\pi _x(f)\w _x$, above, could be
replaced by $f(x)\w _x$, if desired.  Thus, $\Pi (f)$ is the block diagonal
operator, acting on each $V_x$ as the scalar multiplication by $f(x)$.

Also, for each $g$ in $G$, consider the linear operator
  $U_g$
  on $\prodx V_x$,
  given by
  $$
  U_g(\eta )_x = \bool {x\in X_g}\mu _g\big (\eta _{\theta _{g\inv }(x)}\big ),
  \for \eta = (\eta _x)_{x\in X} \in \prodx V_x.
  \equationmark DefineBigUg
  $$

  The above occurence of $\mu _g$ should have actually been written as $\mu
^{\theta _{g\inv }(x)}_g$, but due to the awkward nature of this notation we
will rely on the context to determine the missing superscript.

In what amounts to be essentially a rewording of \ref {MuOnEg}, we have:

\state Proposition \label Functoriality
  Identifying $V_x$ as a subspace of\/ $\prodx V_x$, in the natural way, we
have:
  \izitem
  \zitem if $x\notin X_{g\inv }$, then $U_g$ vanishes on $V_x$,
  \zitem if $x\in X_{g\inv }$, then $U_g$ coincides with $\mu ^x_g$, and hence
maps $V_x$ into $V_{\theta _g(x)}$,
  \zitem if $x\in X_{g\inv }$, then $U_g$ maps $V_x$ bijectively onto $V_{\theta
_g(x)}$,
  \zitem if $x\in X_{g\inv }\cap X_{g\inv h\inv }$, then the composition
  $$
  V_x \flecha {U_g} V_{\theta _g(x)} \flecha {U_h} V_{\theta _{hg}(x)}
  $$
  coincides with $U_{hg}$ on $V_x$.

\Proof (i) and (ii) follow easily by inspection, while (iv) follows directly
from \ref {MuOnEg}.

In order to prove that $U_g$ is bijective from $V_x$ to $V_{\theta _g(x)}$, it
is enough to observe that, by (iv), the restriction of $U_{g\inv }$ to
$V_{\theta _g(x)}$ is the inverse of $U_g$.  \endProof

\state Proposition For every $g$ in $G$, and every $f\in \Lc $, one has that
  $$
  U_g\Pi (f) = \Pi \big (\abar gf\big )U_g.
  $$

\Proof
  Given $\eta = (\w _x)_{x\in X} \in V$, one has for every $x$ in $X$ that
  $$
  \big (U_g\Pi (f) \w \big )_x =
  \bool {x\in X_g}\mu _g\Big (\big (\Pi (f)\w \big ) _{\theta _{g\inv }(x)}\Big
) =
  \bool {x\in X_g}\mu _g\Big (f\big ({\theta _{g\inv }(x)}\big )\eta _{\theta
_{g\inv }(x)}\Big ) \quebra =
  \bool {x\in X_g}f\big (\theta _{g\inv }(x)\big )\mu _g\big (\eta _{\theta
_{g\inv }(x)}\big ) =
  \abar gf|_x(U_g\w )_x =
  \Big (\Pi \big (\abar gf\big )U_g\w \Big )_x.
  $$
  This concludes the proof.  \endProof

As a consequence, there exists a representation $\Pi \times U$ of $\alg $ on
$\prodx V_x$, such that
  $$
  (\Pi \times U)(f\Delta _g) = \Pi (f)U_g, \for f\in \lc {X_{g\inv }}.
  $$

\definition \label Discretization
  The representation $\Pi \times U$ above will be referred to as the
\"{discretization} of the initially given representation $\pi $.

\state Proposition \label IntroQ
  The mapping
  $$
  Q: \xi \in V \mapsto \big (q_x(\xi )\big )_{x\in X} \in \prodx V_x,
  $$
  is injective and covariant relative to the corresponding representations of
$\alg $ on $V$ and on\/ $\prodx V_x$, respectively.

\Proof
  Let $g\in G$, and $f\in \lc {X_g}$.  Then, for every $\xi $ in $V$, and every
$x\in X$, we have
  $$
  \big ((\Pi \times U)(f\d _g)Q(\xi )\big )_x =
  \big (\Pi (f)U_gQ(\xi )\big )_x =
  f(x)\big (U_gQ(\xi )\big )_x =
  f(x)\bool {x\in X_g}\mu _g\big (Q(\xi ) _{\theta _{g\inv }(x)}\big ) \quebra =
  f(x)\mu _g\big (q_{\theta _{g\inv }(x)}(\xi )\big ) =
  f(x)q_x(u_g\xi \big ) =
  q_x\big (\pi (f)u_g\xi \big ) =
  Q\big (\pi (f\d _g)\xi \big )_x.
  $$
  This proves that $Q$ is covariant.

  In order to prove that $Q$ is injective, suppose that $Q(\xi )=0$, for a given
$\xi $ in $V$.  We then claim that, for every $x$ in $X$, there exists a
{\klopen } neighborhood $C$ of $x$, such that
  $$
  \pi (1_{C_x})\xi =0.
  \equationmark LocalValueZero
  $$
  To see this, fixing $x$ in $X$, recall that $q_x(\xi )=0$, by hypothesis, so
$\xi $ lies in $Z_x$ and hence we may write
  $$
  \xi = \sum _{i=1}^n \pi (f_i)\xi _i,
  $$
  where the $f_i$ are in $I_x$, and hence vanish on $x$.  From the fact that the
$f_i$ are locally constant, and finitely many, it follows that there exists a
{\klopen } neighborhood $C_x$ of $x$, where all of the $f_i$ vanish.
Consequently $1_{C_x}f_i=0$, so
  $$
  \pi (1_{C_x})\xi =
  \sum _{i=1}^n \pi (1_{C_x}f_i)\xi _i = 0,
  $$
  proving the claim.

Using \ref {LcNonDeg}, or recycling any one of the above decompositions of $\xi
$, let us again write
  $$
  \xi = \sum _{i=1}^n \pi (f_i)\xi _i,
  $$
  with $f_i\in \Lc $, and $\xi _i\in V$.  Let
  $$
  D = \bigcup _{i=1}^n\supp (f_i),
  $$
  so $D$ is a {\klopen } subset of $X$, and we have
  $$
  \xi =
  \sum _{i=1}^n \pi (1_Df_i)\xi _i =
  \pi (1_D)\sum _{i=1}^n \pi (f_i)\xi _i =
  \pi (1_D)\xi .
  \equationmark XiSuppD
  $$

  Regarding the open cover $\{C_x\}_{x\in X}$ of $D$, where the $C_x$ are as in
the first part of this proof, we may find a finite set $\{x_1,\ldots
,x_p\}\subseteq X$, such that $D\subseteq \bigcup _{i=1}^pC_{x_i}$.  Putting
  $$
  E_k=D\cap C_{x_k}\setminus \bigcup _{i=1}^{k-1}C_{x_i}, \for k=1,\ldots ,p,
  $$
  it is easy to see that the $E_k$ are pairwise disjoint {\klopen } sets, whose
union coincides with $D$.  Observing that $E_k\subseteq C_{x_k}$, we then have
  $$
  \xi \={XiSuppD}
  \pi (1_D)\xi =
  \sum _{k=1}^p \pi (1_{E_k})\xi =
  \sum _{k=1}^p \pi (1_{E_k}1_{C_{x_k}})\xi =
  \sum _{k=1}^p \pi (1_{E_k})\pi (1_{C_{x_k}})\xi \={LocalValueZero} 0.
  $$
  This proves that $Q$ is injective.
  \endProof

As an immediate consequence we have

\state Corollary \label FirstNulSpace
  The null space of\/ $\Pi \times U$ is contained in the null space of\/ $\pi $.

\Proof
  By \ref {IntroQ} we see that $\pi $ is equivalent to a subrepresentation of
$\Pi \times U$, so the conclusion follows.  \endProof

From now on we will consider the subspace
  $$
  \sumx V_x \subseteq \prodx V_x,
  $$
  consisting of the vectors with finitely many nonzero coordinates.  It is easy
to see that this subspace is invariant under $\Pi (f)$, for all $f$ in $\Lc $,
as well as under $U_g$, for all $g$ in $G$, consequently it is also invariant
under $\Pi \times U$.

\state Proposition \label SecondNulSpace
  The null space of the representation obtained by restricting $\Pi \times U$
to\/ $\sumx V_x$ coincides with the null space of\/ $\Pi \times U$ itself.

\Proof
  Given $\cCHb \in \alg $, we must show that if $(\Pi \times U)(\cCHb )$
vanishes on $\sumx V_x$, then it vanishes everywhere.  Writing
  $$
  \cCHb =\sum _{g\in G}f_g\d _g,
  $$
  we have for every $\eta =(\eta _x)_{x\in X}$ in $\prodx V_x$, and for every
$x\in X$, that
  $$
  \big ((\Pi \times U)(\cCHb )\eta \big )_x =
  \sum _{g\in G}\big (\Pi (f_g)U_g \eta \big )_x =
  \sum _{g\in G}f_g(x) \bool {x\in X_g}\mu _g\big (\eta _{\theta _{g\inv
}(x)}\big ).
  $$
  From this we see that the $x^\th $ coordinate of $(\Pi \times U)(\cCHb )\eta $
depends only on the coordinates $\eta _y$, for $y$ of the form $y=\theta _{g\inv
}(x)$, where $g$ is such that $f_g\neq 0$, and $x\in X_g$.  What matters to us
is that the set of $y$'s mentioned above is finite, so if $\eta '$ is defined to
have the same $y$-coordinates as $\eta $, for $y$ on the above finite set, and
zero elsewhere, then $\eta '$ lies in $\sumx V_x$, and
  $$
  \big ((\Pi \times U)(\cCHb )\eta \big )_x = \big ((\Pi \times U)(\cCHb )\eta
'\big )_x = 0.
  $$
  Since $\eta $ and $x$ are arbitrary we deduce that $(\Pi \times U)(\cCHb )=0$,
concluding the proof.  \endProof

As we turn out attention to the restriction of $\Pi \times U$ to $\sumx V_x$, it
is useful to analyze some of its important aspects.
  Initially, regarding the space where it acts, we will identify each $V_x$ as a
subspace of $\sumx V_x$, in the usual way.  Thus, given $\xi $ in $V$, we will
think of $q_x(\xi )$ as the element of $\sumx V_x$ whose coordinates all vanish,
except for the $x^\th $ coordinate which takes on the value $q_x(\xi )$.  Once
this is agreed upon, one may easily show that
  $$\def \quad {\ }
  \matrix {
  \Pi (f)q_x(\xi ) &=& \pi _x(f)q_x(\xi ) &=& q_x\big (\pi (f)\xi \big ), \hfill
\cr \pilar {20pt}
  U_g\big (q_x(\xi )\big ) &=& \bool {x\in X_{g\inv }}\mu _g^x\big (q_x(\xi
)\big ) &=& \bool {x\in X_{g\inv }}q_{\theta _g(x)}(u_g\xi ),
  }
  \equationmark FormulasOnQx
  $$
  for all $f\in \Lc $, $g\in G$, $x\in X$, and $\xi \in V$.

Since $\sumx V_x$ is spanned by the union of the $V_x$, each of which is the
range of the corresponding $q_x$, the formulas above determine the action of the
$\Pi (f)$ and of the $U_g$ on the whole space $\sumx V_x$.  Putting them
together, we may give the following concrete description of the restriction of
$\Pi \times U$ to $\sumx V_x$:

\state Proposition \label ExpressionForPixU
  Given $\cCHb = \sum _{g\in G}f_g\d _g$ in $\alg $, one has for all $x$ in $X$
and $\xi $ in $V$, that
  $$
  (\Pi \times U)(\cCHb )q_x(\xi ) =
  \sum _{g\in G} \bool {x\in X_{g\inv }}q_{\theta _g(x)}\big (\pi (f_g)u_g\xi
\big ).
  $$

\Proof
 The proof is now a simple direct computation:
  $$
  (\Pi \times U)(\cCHb )q_x(\xi ) =
  \sum _{g\in G}\Pi (f_g)U_g\big (q_x(\xi )\big ) =
  \sum _{g\in G}\Pi (f_g) \bool {x\in X_{g\inv }}q_{\theta _g(x)}(u_g\xi )
\quebra =
  \sum _{g\in G} \bool {x\in X_{g\inv }}q_{\theta _g(x)}\big (\pi (f_g)u_g\xi
\big ).
  \closeProof
  $$
  \endProof

Let us now use this to describe the \"{matrix entries} of the operator $(\Pi
\times U)(\cCHb )$ acting on $\sumx V_x$.
  By this we mean that, for each $x$ and $y$ in $X$, we want an expression for
the $y^\th $ component of the vector obtained by applying $(\Pi \times U)(\cCHb
)$ to any given vector in $V_x$, say of the form $q_x(\xi )$, where $\xi \in V$.

The answer is of course the $y^\th $ component of the expression given in \ref
{ExpressionForPixU}, which is in turn given by the partial sum corresponding to
the terms for which $\theta _g(x)=y$.   The desired expression for matrix
entries therefore becomes
  $$
  \big ((\Pi \times U)(\cCHb )q_x(\xi )\big )_y =
  \sum \subsub {g\in G}{\theta _g(x)=y} q_{\theta _g(x)}\big (\pi (f_g)u_g\xi
\big ) =
  q_y \Big (\kern -3pt\sum \subsub {g\in G}{\theta _g(x)=y} \pi (f_g)u_g\xi \Big
).
  \equationmark MatrixEntries
  $$

Recall that in \ref {FirstNulSpace} and \ref {SecondNulSpace} we proved the
following relations among the null spaces of $\pi $, $\Pi \times U$, and the
restriction of the latter to $\sumx V_x$:
  $$
  \Ker (\pi ) \supseteq \Ker (\Pi \times U) = \Ker \big (\Pi \times U|_{\oplus
_{x\in X}V_x}\big ).
  \equationmark SoFarKernelInclusions
  $$

We will now show that equality in fact holds throughout.

\state Theorem \label ThirdNulSpace
  The null space of the representation obtained by restricting\/ $\Pi \times U$
to\/ $\sumx V_x$ coincides with the null space of\/ $\pi $.

\Proof An important aspect of \ref {MatrixEntries}, to be used shortly, is that
since $(\Pi \times U)(\cCHb )$ is well defined on each $V_x$, then so is the
right-hand-side in \ref {MatrixEntries}.  Precisely speaking, if $\xi $ and $\xi
'$ are elements of $V$ such that $q_x(\xi ) = q_x(\xi ')$, then
  $$
  q_y \Big (\kern -3pt\sum \subsub {g\in G}{\theta _g(x)=y} \pi (f_g)u_g\xi \Big
) = q_y \Big (\kern -3pt\sum \subsub {g\in G}{\theta _g(x)=y} \pi (f_g)u_g\xi
'\Big ).
  \equationmark WellDef
  $$

By \ref {SoFarKernelInclusions}, in order to prove the statement, it suffices to
prove that if $\cCHb $ is in the null space of $\pi $, then $(\Pi \times
U)(\cCHb )$ vanishes on $\sumx V_x$, which is the same as saying that its matrix
entries given by \ref {MatrixEntries} vanish for all $x$ and $y$ in $X$.

Again writing $\cCHb =\sum _{g\in G}f_g\d _g$, let $\Gamma $ be the subset of
$G$ consisting of those $g$ for which $f_g\neq 0$, and notice that $\Gamma $
decomposes as the disjoint union of the following subsets:
  $$
  \def \quad {\ }\def \crr {\hfill \stake {12pt}\cr }
  \matrix {
  \Gamma _1 &=& \{g\in \Gamma : y\notin X_g\}, \crr
  \Gamma _2 &=& \{g\in \Gamma : y\in X_g,\ \theta _{g\inv }(y)\neq x\}, \crr
  \Gamma _3 &=& \{g\in \Gamma : y\in X_g,\ \theta _{g\inv }(y)=x\}.\crr
  }
  $$

From our hypothesis that $\pi (\cCHb )=0$, we conclude that, for every $\eta $
in $V$, one has
  $$
  0 = \pi (\cCHb ) \eta =
  \sum _{g\in \Gamma } \pi (f_g\d _g)\eta =
  \sum _{g\in \Gamma } \pi (f_g)u_g\eta .
  \equationmark HypothesisPiVanishes
  $$

This looks enticingly like the last part of \ref {MatrixEntries}, except of
course that here we are summing over all of $\Gamma $, while only the terms
corresponding to $\Gamma _3$ are being considered there.  In order to fix this
discrepancy, notice that $x$ is not a member of the finite set $\{\theta _{g\inv
}(y):g\in \Gamma _2\}$, so we may choose some $\varphi $ in $\Lc $ such that
$\varphi (x)=1$, and $\varphi \big (\theta _{g\inv }(y)\big )=0$, for all $g\in
\Gamma _2$.
  Observing that
  $$
  q_x\big (\pi (\varphi )\xi \big ) \={SummarizePix}
  \varphi (x)q_x(\xi \big ) = q_x(\xi \big ),
  $$
  we will later use \ref {WellDef} in order to replace $\xi $ by
  $$
  \xi ':=\pi (\varphi )\xi
  $$
  in \ref {MatrixEntries}.
  Meanwhile we claim that
  $$
  q_y \big (\pi (f_g)u_g\xi '\big )=0, \for g\in \Gamma _1\cup \Gamma _2.
  \equationmark MeanwhileClaim
  $$
  In order to prove this, observe that
  $$
  q_y \big (\pi (f_g)u_g\xi '\big ) =
  q_y \big (\pi (f_g)u_g\pi (\varphi )\xi \big ) =
  q_y \big (\pi (f_g\abar g\varphi )u_g\xi \big ) \={SummarizePix}
  f_g(y)\abar g\varphi |_y\,q_y (u_g\xi ).
  $$

  If $g\in \Gamma _1$, then the fact that $f_g$ is supported on $X_g$ implies
that $f_g(y)=0$, so the above expression vanishes. On the other hand, if $g\in
\Gamma _2$, then
  $$
  \abar g\varphi |_y = \varphi \big (\theta _{g\inv }(y)\big ) = 0,
  $$
  so the above expression again vanishes, and \ref {MeanwhileClaim} is proved.
Combining this with \ref {HypothesisPiVanishes} we then have
  $$
  0 =
  q_y\Big ( \sum _{g\in \Gamma } \pi (f_g)u_g\xi '\Big ) =
  q_y\Big ( \sum _{g\in \Gamma _1} \pi (f_g)u_g\xi '\Big ) + q_y\Big ( \sum
_{g\in \Gamma _2} \pi (f_g)u_g\xi '\Big ) + q_y\Big ( \sum _{g\in \Gamma _3} \pi
(f_g)u_g\xi '\Big ) \quebra =
  q_y\Big ( \sum _{g\in \Gamma _3} \pi (f_g)u_g\xi '\Big ) \= {WellDef}
  q_y\Big ( \sum _{g\in \Gamma _3} \pi (f_g)u_g\xi \Big ) \={MatrixEntries}
  \big ((\Pi \times U)(\cCHb )q_x(\xi )\big )_y.
  $$
  This shows that $(\Pi \times U)(\cCHb )$ vanishes on $\sumx V_x$, and hence
the proof is concluded.  \endProof

This result will have important consequences for our study of ideals in $\alg $.
The method we shall adopt will be to start with any ideal $J\ideal \alg $, and
then use \ref {IdealIsKernel} and \ref {LocalUnits} to find a representation
$\pi $, as above, such that $\Ker (\pi )=J$.  By \ref {ThirdNulSpace} we may
replace $\pi $ by $\Pi \times U$ acting on $\sumx V_x$, without affecting null
spaces, and it will turn out that the latter is easy enough to understand since
it decomposes as a direct sum of very straightforward sub-representations, which
we will now describe.

\state Proposition \label OrbitInvariant
  Given any $x_0$ in $X$, one has that
  $$
  \bigoplus _{x\in \Orb (x_0)}V_x
  $$
  is invariant under\/ $\Pi \times U$.

\Proof By \ref {Functoriality.ii}, this space is invariant under every $U_g$.
It is also invariant under every $\Pi (f)$, since in fact each $V_x$ has this
property.  Invariance under $\Pi \times U$ then follows.  \endProof

We shall now study the representation obtained by restricting $\Pi \times U$ to
the invariant space mentioned above, so we better give it a name:

\definition \label DefineRepOrbit
  Given $x_0$ in $X$, we shall denote the invariant subspace referred to in \ref
{OrbitInvariant} by
  $
  W_{x_0},
  $
  while the representation of $\alg $ obtained by restricting $\Pi \times U$ to
$W_{x_0}$ will be denoted by $\rho _{x_0}$.

If $R\subseteq X$ is a system of representatives for the orbit relation in $X$,
namely if $R$ contains exactly one point of each orbit relative to the action of
$G$ on $X$, then surely one has
  $$
  \sumx V_x = \bigoplus _{x_0\in R}W_{x_0},
  $$
  while the restriction of $\Pi \times U$ to $\sumx V_x$ is equivalent to
$\bigoplus _{x_0\in R}\rho _{x_0}$.

Before we state the main result of this section we should recall that right
after the proof of \ref {LocalUnits} we fixed an arbitrary ideal $J\ideal \alg
$, which incidentally has been forgotten ever since.

\state Theorem \label JIsIntersRho
  Let $J$ be an arbitrary ideal of $\alg $, and let $\pi $ be a non-degenerate
representation of $\alg $, such that $J=\Ker (\pi )$. Considering the
representations $\rho _x$ constructed above, we have
  $$
  J = \bigcap _{x\in R} \Ker (\rho _x),
  $$
  where $R\subseteq X$ is any system of representatives for the orbit relation
in $X$.

\Proof The null space of $\pi $ coincides with the null space of the restriction
of $\Pi \times U$ to $\sumx V_x$ by \ref {ThirdNulSpace}.  Since the latter
representation is equivalent to the direct sum of the $\rho _x$, as seen above,
the conclusion is evident.  \endProof

\section The representations $\rho _{x_0}$

\label RhoSect In this section we shall keep the setup of the previous section,
such as the ingredients listed in \ref {StandingHyp}, the ideal $J\ideal \alg $,
and the representation $\pi :\alg \to \lin (V)$ fixed there.

The usefulness of Theorem \ref {JIsIntersRho} in describing $J$ is obviously
proportional to the extent to which we may describe the ideals $\Ker (\rho
_{x_0})$ mentioned there, and the good news is that the representations $\rho
_{x_0}$ are well known to us.  In fact they are induced from representations of
isotropy group algebras.
  The main goal of this section is to prove that this is indeed the case.

Our next result refers to the behaviour of the operators
  $$
  U_g:\sumx V_x \to \sumx V_x,
  $$
  when $g$ lies in an isotropy group.

\state Proposition \label EnterKhModule
  Fixing $x_0$ in $X$, let $\Hx $ be the isotropy group of $x_0$. Then, for each
$h$ in $\Hx $, one has that $V_{x_0}$ is invariant under $U_h$.  Moreover, the
restriction of\/ $U_h$ to $V_{x_0}$ is an invertible operator and the
correspondence
  $$
  h\in \Hx \mapsto U_h|_{V_{x_0}}\in \gl (V_{x_0})
  $$
  is a group representation.

\Proof Follows immediately from \ref {Functoriality}.  \endProof

The representation of $\Hx $ on $V_{x_0}$ referred to in the above Proposition
may be integrated to a representation of $\GpAlg \Hx $, which in turn makes
$V_{x_0}$ into a left $\GpAlg \Hx $-module.  Applying the machinery of Section
\ref {InductionSect}, we may then form the induced module $\M \otimes V_{x_0}$,
as in \ref {DefinEinducedModule}, which we may also view as a representation of
$\alg $ on $\M \otimes V_{x_0}$.

\state Theorem \label RhoIsInduced
  For each $x_0$ in $X$, one has that $\rho _{x_0}$ is equivalent to the
representation induced from the left $\KH $-module $V_{x_0}$, as described
above.

\Proof
  Recalling from \ref {DefineRepOrbit} that the space of $\rho _{x_0}$ is
  $$
  W_{x_0} = \bigoplus _{x\in \Orb (x_0)}V_x,
  $$
  consider the bilinear map $T:\M \times V_{x_0}\to W_{x_0}$ given by
  $$
  T\Big (\sum _{k\in \Sx } c_k\delta _k,\xi \Big )=\sum _{k\in \Sx } c_kU_k(\xi
).
  $$
  Recalling that $\M $ is a right $\KH $-module, and viewing $V_{x_0}$ as a left
$\KH $-module via the representation mentioned in \ref {EnterKhModule}, we claim
that $T$ is balanced.  In fact, for every $k\in \Sx $, $h\in \Hx $, and $\xi $
in $V_{x_0}$, one has
  $$
  T(\delta _k\delta _h,\xi )=
  T(\delta _{kh},\xi )=
  U_{kh}(\xi ) \={Functoriality} U_k\big (U_h(\xi )\big ) = T\big (\delta
_k,U_h(\xi )\big ) = T(\delta _k,\delta _h\cdot \xi ).
  $$
  Therefore there exists a unique linear map $\tau :\M \otimes V_{x_0}\to
W_{x_0}$, such that $\tau (\delta _k\otimes \xi )=U_k(\xi )$.  We shall next
prove that $\tau $ is an isomorphism by exhibiting an inverse for it.

With this goal in mind, let $R$ be a system of representatives of left classes
for $\Sx $ modulo $\Hx $.  Thus, if $x$ is in the orbit of $x_0$, there exists a
unique $r$ in $R$ such that $\theta _r(x_0)=x$, so that
  $U_{r\inv }$ maps $V_x$ onto $V_{x_0}$, by \ref {Functoriality}.   We
therefore let
  $$
  \sigma _x:V_x\to \M \otimes V_{x_0}
  $$
  be given by $\sigma _x(\xi ) = \delta _r\otimes U_{r\inv }(\xi )$, for every
$\xi $ in $V_x$.  Putting all of the $\sigma _x$ together, let
  $$
  \sigma :W_{x_0} = \kern -10pt \bigoplus _{x\in \Orb (x_0)}\kern -10pt V_x
\longrightarrow \M \otimes V_{x_0}
  $$
  be the only linear map coinciding with $\sigma _x$ on $V_x$, for every $x$ in
$\Orb (x_0)$.

We claim that $\sigma $ is the inverse of $\tau $.  To see this, let $k$ be any
element of $\Sx $, and let $\xi $ be picked at random in $V_{x_0}$.  Writing
$k=rh$, with $r\in R$, and $h\in \Hx $, set $x=\theta _k(x_0)=\theta _r(x_0)$,
so $U_k(\xi )\in V_x$.  We then have
  $$
  \sigma \big (\tau (\delta _k\otimes \xi )\big ) = \sigma \big (U_k(\xi )\big )
= \delta _r\otimes U_{r\inv }\big (U_k(\xi )\big ) = \delta _r\otimes U_h\xi =
\delta _{rh}\otimes \xi = \delta _k\otimes \xi .
  $$
  This proves that $\sigma \tau $ is the identity on $\M \otimes V_{x_0}$.  On
the other hand, given any $x$ in $\Orb (x_0)$, and any $\xi \in V_x$, write
$x=\theta _r(x_0)$, with $r\in R$, and notice that
  $$
  \tau \big (\sigma (\xi )\big ) = \tau \big (\delta _r\otimes U_{r\inv }(\xi
)\big ) = U_r\big (U_{r\inv }(\xi )\big ) = \xi ,
  $$
  so we see that $\tau \sigma $ is the identity on $W_{x_0}$.

Therefore $\tau $ is an isomorphism between the $\Fld $-vector spaces $\M
\otimes V_{x_0}$ and $W_{x_0}$.  We will next prove that $\tau $ is covariant
for the respective actions of $\alg $, which amount to saying that it is linear
as a map between left ($\alg $)-modules.
  For this, given $g\in G$, and $f\in \lc {X_g}$, we must prove that
  $$
  \tau \big ((f\d _g)\delta _k\otimes \xi \big ) = \rho (f\d _g)\big (\tau
(\delta _k\otimes \xi )\big ), \for k\in \Sx , \for \xi \in V_{x_0}.
  \equationmark ThisIsCovariance
  $$

Given $k$ and $\xi $ as indicated above, the left-hand-side equals
  $$
  \tau \big ((f\d _g)\delta _k\otimes \xi \big ) =
  \bool {gk\in \Sx } f\big (\theta _{gk}(x_0)\big ) \tau (\delta _{gk}\otimes
\xi ) =
  \bool {gk\in \Sx } f\big (\theta _{gk}(x_0)\big ) U_{gk}(\xi ),
  $$
  while the right-hand-side becomes
  $$
  \rho (f\d _g)\big (\tau (\delta _k\otimes \xi )\big ) =
  \Pi (f)U_gU_k(\xi ) = \cdots
  \equationmark RhofDgOnTau
  $$

  Observe that $U_k(\xi )$ is in $V_{\theta _k(x_0)}$, and recall from \ref
{Functoriality} that $U_g$ vanishes on $V_{\theta _k(x_0)}$, unless $\theta
_k(x_0)\in X_{g\inv }$, in which case $U_gU_k$ coincides with $U_{gk}$ on
$V_{x_0}$.  So
  $$
  U_gU_k(\xi ) = \bool {\theta _k(x_0)\in X_{g\inv }}U_{gk}(\xi ).
  $$ Also notice that
  $$
  \theta _k(x_0)\in X_{g\inv } \iff
  \theta _k(x_0)\in X_{g\inv }\cap X_k \iff
  x_0\in \theta _{k\inv }(X_{g\inv }\cap X_k) = X_{k\inv g\inv }\cap X_{k\inv }
\quebra \iff
  x_0\in X_{k\inv g\inv } \iff
  gk\in \Sx ,
  $$
  where we are taking into account that $x_0\in X_{k\inv }$ by default.
 It follows that the expression in \ref {RhofDgOnTau} equals
  $$
  \cdots = \bool {gk\in \Sx } \Pi (f)U_{gk}(\xi ) =
  \bool {gk\in \Sx } f\big (\theta _{gk}(x_0)\big ) U_{gk}(\xi ),
  $$
  because, in the nonzero case, one has that $U_{gk}(\xi )$ lies in $V_{\theta
_{gk}(x_0)}$, and $\Pi (f)$ acts there by scalar multiplication by $f\big
(\theta _{gk}(x_0)\big )$, according to \ref {EigenVec}.  This proves \ref
{ThisIsCovariance}, so $\tau $ is indeed covariant.  \endProof

Summarizing much that we have done so far, the following is the main result of
this work:

\state Theorem \label MainResult
  Let $\theta =(\{\theta _g\}_{g\in G},\{X_g\}_{g\in G})$ be a partial action of
a discrete group $G$ on a Hausdorff, locally compact, totally disconnected
topological space $X$, such that $X_g$ is clopen for every $g$ in $G$.
  Then, every ideal $J\ideal \alg $ is the intersection of ideals induced from
isotropy groups.

\Proof
  Let $R\subseteq X$ be a system of representatives for the orbit relation on
$X$.
  Using \ref {JIsIntersRho} we may write $J$ as the intersection of the null
spaces of the $\rho _x$, for $x$ in $R$,
  while \ref {RhoIsInduced} tells us that $\rho _x$ is equivalent to the
representation induced from a representation of the isotropy group at $x$.  The
null space of $\rho _x$ is therefore induced from an ideal in the group algebra
of said isotropy group by \ref {AnnInducedVsIndIdeal}, whence the result.
  \endProof

Should one want to explicitly write a given ideal $J\ideal \alg $ as the
intersection of induced ideals, the next result should come in handy:

  \def \II {I'} \def \I {I}

\state Proposition
  Under the assumptions of \ref {MainResult}, choose a system $R$ of
representatives for the orbit relation on $X$.
  For each $x$ in $R$, let $\Hx _x$ be the isotropy group at $x$, and let
  $$
  F_x:\alg \to \GpAlg \Hx _x
  $$
  be as in \ref {IntroF}.
  Then, given any ideal $J\ideal \alg $, one has that
  $\II _x:=F_x(J)$ is an admissible ideal of $\GpAlg \Hx _x$, and
  $$
  J = \medcap {x\in R} \Jx {\II _x}.
  $$

\Proof
  That each $\II _x$ is an admissible ideal follows at once from \ref
{FJIsAdmissible}.

  For each $x$ in $R$, let $\I _x$ be the null space of the representation $\rho
_x$ referred to in the proof of \ref {MainResult}, so that
  $$
  J = \medcap {x\in R} \Jx {I_x}.
  $$
  Observe that for each $x\in R$, one has
  $$
  \II _{x} = F_{x}(J) = F_{x}\Big (\,\medcap {y\in R} \Jx {\I _y}\Big )
\subseteq F_{x}(\Jx {\I _{x}}) \explain {IntroIPrime}\subseteq \I _{x}.
  $$
  Consequently
  $\Jx {\II _{x}} \subseteq \Jx {\I _{x}}$, whence
  $$
  \medcap {x\in R} \Jx {\II _x} \subseteq
  \medcap {x\in R} \Jx {\I _x} =J.
  $$

On the other hand, one has by \ref {LargestIdeal} that $\Jx {\II _x}$ is the
largest among the ideals of $\alg $ mapping into $\II _x$ under $F_x$.  Since
$F_x(J)=\II _x$, by definition, we have that $J$ is among such ideals, so
$J\subseteq \Jx {\II _x}$, and then
  $$
  J\subseteq \medcap {x\in R} \Jx {\II _x},
  $$
  concluding the proof.  \endProof

\def \ou {\ \vee \ } \def \implica {\ \imply \ } \def \mirr {meet-irreducib}

\section Primitive, prime and \mirr le ideals

Recall that an ideal $J$ in an algebra $A$ is said to be \"{primitive} if it
coincides with the annihilator of some irreducible module.  It is called
\"{prime} if, whenever $K$ and $L$ are ideals in $A$, then
  $$
  KL\subseteq J \implica (K\subseteq J) \ou (L\subseteq J).
  $$
  Finally, $J$ is said to be \"{\mirr le} if, for any ideals $K$ and $L$ in $A$,
one has
  $$
  K\cap L\subseteq J \implica (K\subseteq J) \ou (L\subseteq J).
  $$

It is well known that every primitive ideal is prime, and since the inclusion
``$KL\subseteq K\cap L$'' holds for any ideals $K$ and $L$, it is clear that
every
  prime ideal
  is
  \mirr le.

The main goal of this section is to show that the induction process preserves
all of the properties mentioned above.

As usual we continue working under \ref {StandingHyp}.

\state Proposition If $I$ is a primitive ideal of\/ $\KH $, then $\JxI $ is
primitive.

\Proof
  By hypothesis $I$ is the annihilator of some irreducible $\KH $-module $V$.
Employing \ref {AnnInducedVsIndIdeal}
  we then have that $\JxI $ is the annihilator of the induced module $\M \otimes
V$, which is irreducible by \ref {IndPreserveIrred}.  Thus $\JxI $ is primitive.
\endProof

In order to deal with primeness and \mirr ility, we first need to prove a
technical result:

\state Lemma \label FIntersection
  Let $J$ and $K$ be ideals in $\alg $.  Then
  \izitem
  \zitem $F(J)\cap F(K)\subseteq F(J\cap K)$.
  \zitem $F(J)F(K)\subseteq F(JK)$.

\Proof We begin by proving (i).  For this, let
  $$
  c = \sum _{h\in \Gamma }c_h\delta _h\in F(J)\cap F(K),
  $$
  where $\Gamma $ is a finite subset of $\Hx $.
  Applying \ref {LiftC} twice, we obtain {\klopen } neighborhoods $V$ and $W$ of
$x_0$, such that
  $
  V,W\subseteq X_h,
  $
  for every $h\in \Gamma $, satisfying
  $$
  c_V:= \sum _{\Gamma }c_h1_V\d _h \in J \and
  c_W:= \sum _{\Gamma }c_h1_W\d _h \in K.
  $$

  Setting $Z=V\cap W$, we have that $Z$ is another {\klopen } neighborhood of
$x_0$, and
  $$
  J\ni 1_Zc_V = \sum _{\Gamma }c_h1_Z1_V\d _h = \sum _{\Gamma }c_h1_Z\d _h
=:c_Z.
  $$
  A similar reasoning shows that $c_Z$ also lies in $K$, so $c_Z\in J\cap K$.
Therefore
  $$
  c = F(c_Z)\in F(J\cap K).
  $$

In order to prove (ii), let $b\in F(J)$ and $c\in F(K)$, and write
  $$
  b=\sum _{h\in \Gamma }b_h\delta _h,\and c=\sum _{h\in \Gamma }c_h\delta _h,
  $$
  where $\Gamma $ is a finite subset of $\Hx $.
  By \ref {LiftC}, there are {\klopen } sets $V$ and $W$, such that
  $
  x_0\in V,W\subseteq X_h,
  $
  for every $h\in \Gamma $, satisfying
  $$
  b_V:= \sum _{\Gamma }b_h1_V\d _h \in J \and
  c_W:= \sum _{\Gamma }c_h1_W\d _h \in K.
  $$
  Observing that $b_V$ and $c_W$ lie in $\cp \Hx $, we then have that
  $$
  bc=F(b_V)F(c_W) =
  \nu \big (E(b_V)\big )\nu \big (E(c_W)\big ) =
  \nu (b_V)\nu (c_W) =
  \nu (b_Vc_W) \quebra =
  \nu \big (E(b_Vc_W)\big ) = F(b_Vc_W) \in F(JK).
  \closeProof
  $$ \endProof

We may now prove the result announced earlier:

\state Theorem Let $I$ be an ideal in $\KH $.  If $I$ is prime or \mirr le, then
so is $\JxI $.

\Proof
  Let us first address \mirr ility, so
  suppose that $K$ and $L$ are ideals in $\alg $, such that $K\cap L\subseteq
\JxI $.  Then
  $$
  F(K)\cap F(L)\explain {FIntersection}\subseteq
  F(K\cap L)\subseteq F(\JxI ) \explain {IntroIPrime.ii}\subseteq I.
  $$

Assuming that $I$ is \mirr le, we have that either $F(K)$ or $F(L)$ is contained
in $I$.  Supposing without loss of generality that the first alternative is
true, that is, $F(K)\subseteq I$, we then have
  $$
  K\explain {LargestIdeal.i}\subseteq \Jx {F(K)}\subseteq \JxI ,
  $$
  so $\JxI $ is \mirr le.  The proof of the result for prime ideals is obtained
by going through the present proof, replacing all intersections with products.
\endProof

\section Topologically free points

\label TopFreeSect As we already hinted upon, topologically free minimal actions
prevent the appearance if nontrivial induced ideals.  In this section we wish to
further explore this aspect.  We keep enforcing \ref {StandingHyp}.

\definition \label DefTopFree
  \izitem
  \zitem We say that $\theta $ is a \"{topologically free} partial action if,
for every $g$ in $G\setminus \{1\}$, the fixed point set
  $$
  F_g := \{x\in X_{g\inv }: \theta _g(x)=x\}
  $$
  has empty interior.
  \zitem We shall say that a point $x_0$ in $X$ is \"{topologically free} if,
for every $g$ in $G\setminus \{1\}$, and every open set $V$, with $x_0\in
V\subseteq X_{g\inv }$, there exists some $y\in V\cap \Orb (x_0)$, such that
$\theta _g(y)\neq y$.

If $x_0$ is not fixed by $\theta _g$, then the point $y$ referred to in \ref
{DefTopFree.ii} may clearly be taken to be $x_0$ itself, so the condition is
automatically satisfied for such a $g$.  In other words, this condition is only
relevant for $g$ in the isotropy group of $x_0$.

Another way to describe the notion of topologically free point is to say that
there is no subset of $\Orb (x_0)$ containing $x_0$, open in the relative
topology, and consisting of fixed points for a nontrivial group element $g$.

Given the relative notion of the concept of ``interior'', one may find a
topologically free partial action admiting a invariant subspace $Y\subseteq X$,
such that the restriction of $\theta $ to $Y$ is no longer topologically free.
However it is clear that the notion of topologically free \"{point} is not
affected by restricting the action to an invariant subset, as long as the point
under consideration lies in such a subset.

Still another equivalent description of topologically free points is given by
the following:

\state Proposition \label TopFreeOnOrb
  Given $x_0$ in $X$, the following are equivalent:
  \izitem
  \zitem $x_0$ is topologically free,
  \zitem the restriction of\/ $\theta $ to $\cOrb (x_0)$ (the closure of the
orbit of $x_0$) is a topologically free
  partial action.

\Proof (i) $\imply $ (ii).
  Since we are only concerned with $\cOrb (x_0)$, rather than the whole of $X$,
we may replace the latter by the former, and hence assume that the orbit of
$x_0$ is dense in $X$.  As already observed, this restriction does not affect
condition (i).

  Assume by contradiction that $g$ is a nontrivial group element whose fixed
point set
  $
  F_g
  $
  has a nontrivial interior, so there exists a nonempty open set $V\subseteq
F_g$.  Since the orbit of $x_0$ is assumed to be dense, there is some $k$ in
$\Sx $ such that $\theta _k(x_0)\in V$.   It is then easy to prove that $\theta
_{k\inv gk}$ is the identity on the open set
  $$
  U:=\theta _{k\inv }(X_k\cap V),
  $$
  which contains $x_0$.
  In particular $\theta _{k\inv gk}$ is the identity on $U\cap \Orb (x_0)$,
hence contradicting (i).

\medskip \noindent (ii) $\imply $ (i).
  Again by contradiction, assume that $1\neq g\in G$, and that $V$ is an open
set with $x_0\in V\subseteq X_{g\inv }$, and $V\cap \Orb (x_0)\subseteq F_g$.
We then claim that
  $$
  V\cap \cOrb (x_0)\subseteq F_g,
  $$
  as well.  To see this, let $y\in V\cap \cOrb (x_0)$. Then $y$ is the limit of
a net $\{u_i\}_i\subseteq \Orb (x_0)$, and since $y\in V$, we have that $u_i\in
V$ for all sufficiently large $i$.  For such $i$'s we have $u_i\in V\cap \Orb
(x_0)\subseteq F_g$, so
  $$
  \theta _g(y)=\lim _i \theta _g(u_i)=\lim _i u_i = y,
  $$
  proving that $y\in F_g$.  The claim is therefore proven, contradicting (ii).
\endProof

As a consequence we see that two points in $X$ having the same orbit closure are
either both topologically free or both fail to satisfy this property.

Topologically free points actually enjoy a slightly stronger property as
described next:

\state Proposition \label StrongerProperty
  Let $x_0$ be a topologically free point, let $\Gamma $ be a finite subset of\/
$G\setminus \{1\}$, and let $V$ be an open set with
  $$
  x_0\in V\subseteq \raise 3pt \hbox {$\ds \bigcap _{g\in \Gamma }$}X_{g\inv }.
  $$
  Then
  there exists some $y\in V\cap \Orb (x_0)$, such that $\theta _g(y)\neq y$, for
all $g$ in $\Gamma $.

\Proof
  By restricting $\theta $ to the closure of the orbit of $x_0$ we may assume
that $\Orb (x_0)$ is dense in $X$.

For each $g$ in $\Gamma $, let
  $$
  \Phi _g=F_g\cap V = \{x\in V:\theta _g(x)=x\}.
  $$
  Then clearly $\Phi _g$ is a closed subset (relative to $V$) and by \ref
{TopFreeOnOrb} we have that $\Phi _g$ has no interior (relative to $X$, and
hence also relative to $V$).  Consequently $\bigcup _{g\in G}\Phi _g$ is a
closed set with empty
  interior\fn {A finite union of closed sets with empty interior always has
empty interior.},
  whence
  $$
  V\setminus \bigcup _{g\in G}\Phi _g
  $$
  is a nonempty open set (relative to $V$ and hence also relative to $X$).
Since we are assuming that the orbit of $x_0$ is dense, we conclude that there
is some $y$ in said orbit which also lies in the above open set.  This concludes
the proof.
  \endProof

The conflict between topological freeness and induced ideals is clearly
expressed by the following:

\state Proposition \label TopFreeNoAdmiss
  When $x_0$ is topologically free, the only admissible ideals of\/ $\KH $ are
the trivial ones, namely
 $\{0\}$ and $\KH $, itself.  Consequently the only induced ideals arising from
ideals in $\KH $ are the trivial ones described in \ref {Examples}.

\Proof
  Let $I\ideal \KH $ be a nonzero admissible ideal.  We first claim that there
exists some
  $
  c=\sum _{h\in H}c_h\delta _h\in I,
  $
  with $c_1\neq 0$.  To see this let
  $d=\sum _{h\in H}d_h\delta _h$ be any nonzero element of $I$.  Choose $h_0$ in
$H$ such that $d_{h_0}\neq 0$, and let
  $$
  c=d\delta _{h_0\inv } = \sum _{h\in H}d_h\delta _{hh_0\inv },
  $$
  so that $c$ is also in $I$, and
  $
  c_1=d_{h_0}\neq 0,
  $
  proving the claim.

Working with $c$, choose $V$ as in \ref {AdmissibleCharact}.  Letting
  $$
  \Gamma =\{h\in H: c_h\neq 0\},
  $$
  we may clearly assume that $V\subseteq \bigcap _{h\in \Gamma }X_{h\inv }$.
Employing \ref {StrongerProperty}, let $y$ be an element of the orbit of $x_0$,
belonging to $V$, and not fixed by any $h\in \Gamma \setminus \{1\}$.

Writing $y=\theta _k(x_0)$, we claim that $\Gamma \cap kHk\inv =\{1\}$.  To see
this it is enough to observe that $\theta _k(x_0)$ is fixed by the elements of
$kHk\inv $, while the only element of $\Gamma $ having this property is the
unit.

By \ref {AdmissibleCharact} we then conclude that
  $$
  I \ni
  \delta _{k\inv }\Big (\sum _{h\in \Gamma \cap kHk\inv }c_h\delta _h\Big
)\delta _k =
  \delta _{k\inv }(c_1\delta _1)\delta _k = c_1\delta _1.
  $$
  This implies that $I$ contains a nonzero multiple of the unit $\delta _1$, an
invertible element, whence $I=\KH $, concluding the proof.

Regarding the last sentence in the statement, if $I$ is any ideal in $\KH $,
then by \ref {ChangeIdeals} there exists an admissible ideal $I'$ such that
$\JxI = \Jx {I'}$.  By the first part of the proof we have that $I'$ is either
$\{0\}$ or $\KH $, as desired.  \endProof

\section Regular Points

In this section, still under \ref {StandingHyp}, we will study points possessing
a property which may be considered as being in the other end of the spectrum,
relative to topological freeness.

\definition
  A point $x_0$ in $X$ is said to be
  \"{strongly regular} (resp.~\"{regular}) if, for every $h$ in the isotropy
group of $x_0$, there exists an open set $V$ with $x_0\in V\subseteq X_{h\inv
}$, and such that $\theta _h$ is the identity on $V$ (resp.~on $V\cap \Orb
(x_0)$).

Like the notion of topological free point, the notion of regular point given
above is phrased in such a way as to depend only on the action of $G$ on the
orbit of the point under consideration (seen under the relative topology).  This
has the advantage of being mostly an atribute of the point, rather than of the
action.  However, the same cannot be said of the notion of strongly regular
point.  In any case, it is easy to see that every strongly regular point is also
regular.

The following result is the partial actions version of (and it follows from)
\cite [Lemma 3.3.a]{ReznikofEtAll}.

\state Proposition If $G$ is countable, then the set of strongly regular points
is dense.

\Proof Observe that a point $x_0$ fails to be strongly regular precisely when it
lies in the fixed point set
  $
  F_g
  $
  for some $g$ in $G$, but it does not belong to the interior of $F_g$.  This is
obviously to say that $x_0\in \partial  F_g$, meaning the boundary of $F_g$.  So
the set of points which are not strongly regular is precisely the set
  $$
  {\cal S} := \bigcup _{g\in G}\partial F_g.
  $$

On the other hand, since $F_g$ is closed, its boundary is a closed set with
empty interior.  Therefore, should $G$ be countable, we have that $\cal S$ is of
first category in Baire's sense, hence its complement, namely the set of
strongly regular points, is dense.  \endProof

For regular points, a much simpler characterization of admissibility may be
given, if compared to \ref {AdmissibleCharact}.  This will be done based on a
simpler decoding of the information that ``$c_V\in \JxI $'' in the first
paragraph of the proof of \ref {AdmissibleCharact}.  In order to highlight this
simplification, which will be used elsewhere later, we will isolate the
technicalities involved in the next two auxiliary results.

\state Lemma \label FixingConjuga
  Let $\Gamma $ be any subset of $G$, and let $k\in \Sx $ be such that $\theta
_k(x_0)$ is fixed by $\theta _g$, for all $g$ in $\Gamma $.  Then
  \izitem
  \zitem $\Gamma \subseteq k\Hx k\inv $,
  \zitem for every $l$ in $G$ such that $\Gamma \cap k\Hx l\inv $ is nonempty,
one has that
  $l\in \Sx $, that $\theta _l(x_0)=\theta _k(x_0)$, and moreover
  $\Gamma \subseteq k\Hx l\inv $.

\Proof
  (i) Given $g\in \Gamma $, we have that
  $\theta _g\big (\theta _k(x_0)\big )=\theta _k(x_0)$, so
  $\theta _{k\inv }\big (\theta _g(\theta _k(x_0))\big )=x_0$, whence $k\inv g
k$ is in $\Hx $, and consequently $g \in k\Hx k\inv $.  This proves (i).

\medskip \noindent  (ii) Pick
  $g_1$ in $\Gamma \cap k\Hx l\inv $,
  so that
  $
  h_1:=k\inv g_1 l\in \Hx .
  $
  Therefore
  $$
  \theta _k(x_0) = \theta _{g_1\inv }(\theta _k(x_0)\big ) = \theta _{g_1\inv
}\big (\theta _k(\theta _{h_1}(x_0))\big ) = \theta _{g_1\inv kh_1}(x_0) =
\theta _l(x_0).
  $$
  This proves that $l$ is in $\Sx $, and that $\theta _k(x_0)=\theta _l(x_0)$.
Next, picking any $g\in \Gamma $, notice that
  $$
  x_0 =
  \theta _{k\inv }\big (\theta _k(x_0)\big ) =
  \theta _{k\inv }\big (\theta _g(\theta _k(x_0))\big ) =
  \theta _{k\inv }\big (\theta _g(\theta _l(x_0))\big ) =
  \theta _{k\inv g l}(x_0).
  $$
  Thus $k\inv g l\in \Hx $, and so $g\in k\Hx l\inv $, proving (ii).

\endProof

\state Lemma \label LemaOnMembership
  Let $I$ be an ideal in $\KH $, and let
  $
  c=\sum _{g\in \Gamma }c_g\delta _g,
  $
  be an arbitrary element of\/ $\GpAlg G$, where $\Gamma $ is a finite subset of
$G$. Suppose that $V$ is a {\klopen } set such that
  $V\subseteq X_g$,
  and $\theta _{g\inv }$ coincides with the identity on $V\cap \Orb (x_0)$,
  for all $g$ in $\Gamma $.
  Then
  $$
  c_V:= \sum _{g\in \Gamma }c_g1_V\d _g
  $$
  lies in $\JxI $ if and only if, for every $k$ in $\Sx $, such that $\theta
_k(x_0)\in V$, one has that
  $$
  \delta _{k\inv }c\delta _k \in I.
  $$

\Proof
  We begin with the ``only if'' part, so we assume that $c_V\in \JxI $.  Given
$k$ in $\Sx $, with $\theta _k(x_0)\in V$, we have
  $$
  I\ni \langle \delta _k,c_V\delta _k\rangle \={DeltaKbDeltaL} \sum _{g\in
\Gamma \cap k\Hx k\inv }c_g1_V\big (\theta _k(x_0)\big )\, \delta _{k\inv gk} =
  \delta _{k\inv }\Big (\sum _{g\in \Gamma \cap k\Hx k\inv }c_g \delta _g\Big
)\delta _k.
  $$

  Observing that $\theta _k(x_0)$ lies in $V\cap \Orb (x_0)$, we have by
hypotheses that $\theta _k(x_0)$ is fixed by $\theta _{g\inv }$, and hence also
by $\theta _g$, for every $g$ in $\Gamma $.  We then conclude from \ref
{FixingConjuga.i} that $\Gamma \subseteq k\Hx k\inv $, so the computation above
gives
  $
  \delta _{k\inv }c\delta _k\in I,
  $
  as desired.

In order to prove the ``if'' part, let us show that $c_V$ lies in $\JxI $ by
employing the criteria given in \ref {UseDeltaKbDeltaL}.  For this we must prove
that, for every $k$ and $l$ in $\Sx $, one has that
  $$
  \sum _{g\in \Gamma \cap k\Hx l\inv }c_g 1_V\big (\theta _k(x_0)\big )\delta
_{k\inv gl} \in I.
  \equationmark LemaOnMembershipTwo
  $$

  There are two situations in which the above vanishes, in which case there is
nothing to do, namely when
    $\Gamma \cap k\Hx l\inv $ is the empty set,
    or when
    $\theta _k(x_0)\notin V$.
    Ignoring these, let us assume that the opposite is true, namely that $\Gamma
\cap k\Hx l\inv $ is nonempty and that $\theta _k(x_0)$ lies in $V$, which in
turn implies that $\theta _k(x_0)$ is fixed by $\Gamma $.  Therefore \ref
{FixingConjuga.ii} gives $\theta _k(x_0)=\theta _l(x_0)$, and $\Gamma \subseteq
k\Hx l\inv $. We then see that the term appearing in \ref {LemaOnMembershipTwo}
is given by
  $$
  \sum _{g\in \Gamma }c_g \delta _{k\inv gl} =
  \delta _{k\inv }\Big (\sum _{g\in \Gamma }c_g\delta _g\Big )\delta _l =
  \delta _{k\inv }c\,\delta _l = \delta _{k\inv }c\,\delta _k\delta _{k\inv l}.
  $$

  To see that this lies in $I$, notice that
  $k\inv l\in \Hx $, because $\theta _k(x_0)=\theta _l(x_0)$, and moreover that
$\delta _{k\inv }c\,\delta _k$ is in $I$ by hypothesis.  So \ref
{LemaOnMembershipTwo} follows from the fact that $I$ is an ideal in $\KH $.  We
then conclude that $c_V\in \JxI $, thanks to \ref {UseDeltaKbDeltaL}.  \endProof

The promissed simplified characterization of admissibility is given next:

\state Proposition \label ConditionForRegAdmiss
  Suppose that $x_0$ is regular.  Then an ideal $I\ideal \KH $ is admissible if
and only if, for every $c$ in $I$, there exists a neighborhood $V$ of $x_0$,
such that
  $$
  \delta _{k\inv }c\delta _k \in I,
  $$
  for all $k$ in $\Sx $, such that $\theta _k(x_0)\in V$.

\Proof
  We begin exactly as in the proof of \ref {AdmissibleCharact}:
  supposing that $I$ is admissible, pick $c=\sum _{h\in \Gamma }c_h\delta _h$ in
$I$, where $\Gamma $ is a finite subset of $\Hx $.   By hypothesis $c$ is in
$F\big (\JxI \big )$, so \ref {LiftC} provides a {\klopen } neighborhood $V$ of
$x_0$, such that
  $$
  c_V:= \sum _{h\in \Gamma }c_h1_V\d _h \in \JxI .
  $$

  Given that $x_0$ is regular, and upon shrinking $V$, if necessary, we may
assume that $\theta _{h\inv }$ is the identity on $V\cap \Orb (x_0)$, for every
$h$ in $\Gamma $.  The conclusion then follows from \ref {LemaOnMembership}.

Conversely, assuming that $I$ satisfies the condition in the statement, let us
prove that $I$ is admissible, namely that $F\big (\JxI \big )\supseteq I$.  So,
pick any $c=\sum _{h\in \Gamma }c_h\delta _h$ in $I$, where $\Gamma $ is a
finite subset of $\Hx $.  Using the hypothesis, we then choose $V$ as in the
statement, which we may clearly suppose to be {\klopen }.   Again because $x_0$
is regular, we may assume that $\theta _{h\inv }$ is defined and coincides with
the identity on $V\cap \Orb (x_0)$, for every $h$ in $\Gamma $.

By hypothesis, and by \ref {LemaOnMembership}, if follows that $c_V\in \JxI $,
whence
  $$
  c = F(c_V)\in F\big (\JxI \big ),
  $$
  as desired.  \endProof

An important, albeit trivial conclusion to be drawn from the above result is:

\state Corollary \label CommutativeAdmiss
  If $G$ is commutative and $x_0$ is a regular point of $X$, then every ideal
of\/ $\KH $ is admissible.

\def \ad {\hbox {\sl Ad}} \def \pad {\hbox {\sl pAd}}     \def \H {\Hx } \def \Aut {\hbox {Aut}}

\section Normal ideals

It is interesting to notice that, while the admissibility condition given in
\ref {ConditionForRegAdmiss} is a combination of dynamical features
(viz.~``$\theta _k(x_0)\in V$") and algebraic properties (viz.~``$\delta _{k\inv
}c\delta _k \in I\,$"), the algebraic properties alone ensure admissibility in
\ref {CommutativeAdmiss}.

In this section we shall discuss other purely algebraic conditions on an ideal
of $\KH $ which are enough to guarantee admissibility, regardless of any other
dynamical restrictions.

Given a group $G$ and a field $\Fld $, recall that the well known \"{adjoint}
action of $G$ on $\GpAlg G$ is the map
  $$
  \ad :G\to \Aut (\GpAlg G)
  $$
  given by
  $$
  \ad _g(a)=\delta _ga\delta _{g\inv }, \for g\in G, \for a\in \GpAlg G.
  $$

Given any subgroup $H$ of $G$, observe that $\GpAlg H$ is invariant under $\ad $
if and only if $H$ is a normal subgroup.  Regardless of normality, we may always
restrict $\ad $ to a \"{partial action} of $G$ on $\GpAlg H$, as in \cite
[3.2]{mybook}.  The main ingredients of this construction are as follows: for
each $g$ in $G$, we let
  $$
  D_g=\GpAlg H\cap \ad _g(\GpAlg H),
  $$
  and we let
  $$
  \pad _g:D_{g\inv }\to D_g
  $$
  be the restriction of $\ad _g$ to $D_{g\inv }$.  It is well known that $\pad $
is then a partial action (in the category of sets).

\definition The above partial action will be called the \"{adjoint partial
action} of $G$ on $\GpAlg H$.

It is easy to see that each $D_g$ is a subalgebra of $\GpAlg H$, while the $\pad
_g$ are algebra isomorphisms. However $\pad $ cannot be viewed as an algebraic
partial action, as defined in \cite [6.4]{mybook}, because the $D_g$ are not
ideals in $\GpAlg H$, but alas, $\pad $ is a legitimate set theoretical partial
action cf.~\cite [2.1]{mybook}.

\definition Let $H$ be a subgroup of a group $G$, and let $I$ be an ideal in
$\GpAlg H$.  We shall say that $I$ is \"{normal relative to} $G$, if $I$ is
invariant \cite [2.9]{mybook} under the adjoint partial action of $G$ on $\GpAlg
H$.

Thus, to say that $I$ is normal is to say that for every $c\in I$, and every $g$
in $G$ such that $\delta _gc\delta _{g\inv }\in \GpAlg H$, one has that $\delta
_gc\delta _{g\inv }\in I$.

One should view this as the best possible effort made by the ideal $I$ in trying
to embrace all element of the above form $\delta _gc\delta _{g\inv }$, except of
course that this is impossible in the hopeless cases when such elements are not
even in $\GpAlg H$!

\state Proposition \label NormalIsAdmissible
  Under \ref {StandingHyp}, let $x_0$ be a regular point of $X$, and let $\Hx $
be its isotropy group.  Then every ideal $I\ideal \GpAlg H$ which is normal
relative to $G$, is also admissible.

\Proof
  We will verify the conditions of \ref {ConditionForRegAdmiss}.  Thus, given
$c$ in $I$, write $c=\sum _{h\in \Gamma }c_h\delta _h$, where $\Gamma \subseteq
H$ is a finite set, and choose a neighborhood $V$ of $x_0$, such that $\theta
_h$ is the identity map on $V\cap \Orb (x_0)$, for every $h$ in $\Gamma $.
Still focusing on \ref {ConditionForRegAdmiss}, pick any $k$ in $\Sx $ such that
$\theta _k(x_0)\in V$.

  We then claim that $\delta _{k\inv }c\delta _k\in \KH $.  To see this, notice
that, for every $h$ in $\Gamma $, one has that $\theta _h$ fixes $\theta
_k(x_0)$, meaning that
  $
  \theta _h\big (\theta _k(x_0)\big ) = \theta _k(x_0),
  $
  from where we deduce that
  $$
  \theta _{k\inv }\big (\theta _h(\theta _k(x_0))\big ) = x_0.
  $$
  So $k\inv h k\in \Hx $, whence
  $$
  \delta _{k\inv }c\delta _k = \sum _{h\in \Gamma }c_h\delta _{k\inv h k} \in
\KH .
  $$
  The invariance of $I$ under the adjoint partial action then implies that
$\delta _{k\inv }c\delta _k\in I$, concluding the verification of the conditions
of \ref {ConditionForRegAdmiss}, and hence proving that $I$ is admissible.
\endProof

A source of examples of normal ideals is as follows:

\state Proposition \label SourceOfNormal
  Let $H$ be a subgroup of a group $G$, and let $J$ be any ideal in $\GpAlg G$.
Then the ideal $I$ of\/ $\GpAlg H$ given by
  $$
  I=J\cap \GpAlg H
  $$
  is normal relative to $G$.

\Proof
  If $c$ is in $I$, then for every $g$ in $G$, one has that $\delta _gc\delta
_{g\inv }\in J$.  If the latter happens to also lie in $\GpAlg H$, then it
clearly belongs to $I$.  Therefore $I$ is normal.
  \endProof

A concrete example is the \"{augmentation ideal} $I_H$ given by
  $$
  I_H = \Ker (\varepsilon _H,)
  $$
  where $\varepsilon _H$ is the \"{augmentation map}, namely the map
  $\varepsilon _H:\GpAlg H\to \Fld $, given by
  $$
  \varepsilon _H\Big (\sum _{h\in H} c_h\delta _h\Big ) = \sum _{h\in H} c_h.
  $$

\state Proposition \label AugmAdmiss
  Let $H$ be a subgroup of a group $G$.  Then the augmentation ideal $I_H$ is
normal relative to $G$.

\Proof
  This ideal being the intersection of $\GpAlg H$ with the augmentation ideal
$I_G$ of $G$, the conclusion follows from \ref {SourceOfNormal}.
  \endProof

Incidentally, the ideal referred to in \cite {LisaEtAll} is related to the ideal
induced by $I_H$.  In particular we have:

\state Proposition {\rm (cf.~\cite {LisaEtAll})}
  Assuming \ref {StandingHyp} and that $\alg $ is simple, one has that $\theta $
is topologically free.

\Proof
  Suppose by contradiction that $\theta $ is not topologically free. Then there
exists a nontrivial $g$ in $G$ whose fix point set
  $
  F_g
  $
  has nonempty interior.  Since the regular points are dense in $X$, we may pick
a regular point $x_0$ in $V$.  In particular $g$ lies in the isotropy group $\Hx
$ of $x_0$, so $\Hx $ is a nontrivial group, whence
  $$
  \{0\}\subsetneq I_\Hx \subsetneq \KH ,
  $$
  where $I_\Hx $ is the augmentation ideal of $\KH $.  Observe that the three
ideals above are admissible by \ref {TrivialIdeals}, \ref {AugmAdmiss} and \ref
{NormalIsAdmissible}, so by the uniqueness part of \ref {ChangeIdeals}, we have
  $$
  \Jx {\{0\}}\subsetneq \Jx {I_\Hx }\subsetneq \Jx {\KH }\,.
  $$

However, since $\alg $ is supposed to be a simple algebra, it is impossible to
find three distinct ideals as above.  This is a contradiction, and hence the
statement is proved.
  \endProof

In view of \ref {NormalIsAdmissible} one could ask whether conditions can be
found regarding an ideal $I\ideal \KH $, which would ensure $I$ to be admissible
regardless of any dynamical condition, as in \ref {NormalIsAdmissible}, but also
regardless of $x_0$ being a regular point.
  Except for the trivial ideals treated in \ref {Examples}, this seems to be
impossible in view of \ref {TopFreeNoAdmiss}, where topological freeness, an
eminently dynamical condition, overrides any algebraic condition one could think
of.

\def \pconj {\Psi }

\def \nameDecoration {hat} \def \decoration #1{\csname \nameDecoration
\endcsname {#1}} \def \alt #1#2{\ifnum #2=0 #1\else \decoration {#1}\fi } \def
\I #1{\alt I#1} \def \H #1{\alt \Hx #1} \def \Ga #1{\GpAlg \H #1} \def \F
#1{\alt F#1} \def \S #1{\alt \Sx #1}
  \def \J #1#2{\hbox {\sl I}{\alt {\hbox {\sl n}}#1\hbox {\sl d}(#2})}

\def \x #1{\alt x#1_0} \def \T #1{\alt {{\cal T}}#1}

\section Transposition

\label TranspositionSect So far we have concentrated our study on induced ideals
relative to a single point $\x 0$ in $X$, but now we would like to conduct a
comparative study.  So, besides assuming \ref {StandingHyp}, and hence having
fixed a point $x_0$, we will fix another point in $X$, denoted $\x 1$, and we
will discuss the relationship between ideals induced relative to $\x 0$ and its
peer $\x 1$.

  \def \qt #1#2{``$#1_{\x #2}\kern -1pt$''}Having two points in sight,
  it is now crucial that we distinguish the sets $\Hx $ and $\Sx $ introduced in
\ref {BasicNotations}, depending on whether $\x 0$ or $\x 1$ is concerned.  One
alternative would be to employ their official notation with corresponding
subscripts, such as
  \qt {\Hx }0, \qt {\Sx }0, \qt {\Hx }1 and \qt {\Sx }1.
  However we will really only consider the induction process for the two points
$\x 0$ and $\x 1$ chosen above, so we will prefer to save on notation by keeping
the undecorated notation when $\x 0$ is considered, and writing
  $\H 1$ and $\S 1$, when we are talking about $\x 1$.

The maps $E$, $\nu $ and $F$, respectively introduced in \ref {IntroE}, \ref
{IntroNu} and \ref {IntroF}, also need to be distinguished, so we will adopt the
above policy of decorating everything regarding $\x 1$ with a ``\nameDecoration
''.

Finally, the induction process itself needs to be distinguished, so we will
write $\J 1{\I 1}$, if inducing an ideal $\I 1\ideal \Ga 1$, relative to $\x 1$,
while retaining our previous notation regarding $\x 0$.

The crucial way in which the two induction processes are related may be subsumed
by a correspondence between ideals in $\GpAlg \H 0$ and ideals in $\GpAlg \H 1$,
defined as follows: given an ideal $\I 0\ideal \Ga 0$, we may form the induced
ideal $\J 0{\I 0}$, and then we have by \ref {FJIsAdmissible} that $\F 1\big (\J
0{\I 0}\big )$ is an admissible ideal in $\Ga 1$ (relative to $\x 1$, of
course).

\definition \label DefineTransp
  Given an ideal $\I 0\ideal \Ga 0$, we shall let
  $$
  \T 1(\I 0)=\F 1\big (\J 0{\I 0}\big ),
  $$
  so that $\T 1$ is a map from
  the set of all ideals in $\GpAlg \H 0$ into
  the set of all admissible ideals in $\GpAlg \H 1$.
  We shall refer to $\T 1(\I 0)$ as the \"{transposition of $\I 0$ from $\Ga 0$
to $\Ga 1$}.
  Likewise, given an ideal $\I 1\ideal \Ga 1$, its transposition from $\Ga 1$ to
$\Ga 0$ is defined by
  $$
  \T 0(\I 1)=\F 0\big (\J 1{\I 1}\big ).
  $$

Since we are in the business of studying induced ideals we don't really care so
much about non admissible ideals, so we will shortly restrict ourselves to
transposing admissible ideals only.  Nevertheless one might observe that an
ideal $\I 0\ideal \Ga 0$ is admissible if and only if it coincides with its own
transposition from $\Ga 0$ to itself.

Even before we fully understand the transposition map, we may prove a few
important facts:

\state Proposition \label InclusionInduced
  Let $\I 0\ideal \Ga 0$ and $\I 1\ideal \Ga 1$ be admissible ideals, then the
following are equivalent
  \izitem
  \zitem $\J 0{\I 0}\subseteq \J 1{\I 1}$,
  \zitem $\T 1(\I 0)\subseteq \I 1$.
  \medskip \noindent In addition, when the above equivalent conditions hold, and
both $\I 0$ and $\I 1$ are proper ideals, then $\cOrb (\x 0)\supseteq \cOrb (\x
1)$.

\Proof (i) $\imply $ (ii):
   We have
  $$
  \T 1(\I 0) =
  \F 1\big (\J 0{\I 0}\big ) \subseteq
  \F 1\big (\J 1{\I 1}\big ) =
  \I 1,
  $$
  where the last equality is a consequence of the fact that $\I 1$ is
admissible.

  \bigskip \noindent (ii) $\imply $ (i): Observing that our hypothesis reads $\F
1\big (\J 0{\I 0}\big )\subseteq \I 1$, recall from \ref {LargestIdeal} that $\J
1{\I 1}$ is the largest ideal mapping into $\I 1$ under $\F 1$, whence (i)
holds.

Regarding the last sentence in the statement, we have by \ref {InterCoefAlg}
that the intersection $\JxI \cap \,\Lc $ consists of all $f$ in $\Lc $ vanishing
on $\cOrb (\x 0)$.  Therefore (i) implies that every such $f$ necessarily also
vanishes on $\cOrb (\x 1)$, from where the conclusion follows.  \endProof

The fact that $\T 1(\I 0)=\I 1$ is not equivalent to $\I 0=\T 0(\I 1)$, so our
result for equality of induced ideals must mention both:

\state Theorem \label MainEqualityInduced
  Let $\I 0\ideal \Ga 0$ and $\I 1\ideal \Ga 1$ be admissible ideals, then the
following are equivalent:
  \izitem
  \zitem $\J 0{\I 0}=\J 1{\I 1}$,
  \zitem $\T 1(\I 0)=\I 1$, and $\I 0=\T 0(\I 1)$,
  \zitem $\T 1(\I 0)\subseteq \I 1$, and $\I 0\supseteq \T 0(\I 1)$.

\Proof (i) $\imply $ (ii): We have
  $$
  \T 1(\I 0) = \F 1\big (\J 0{\I 0}\big ) = \F 1\big (\J 1{\I 1}\big )
\={DefAdmiss} \I 1,
  $$
  and one similarly proves that
  $\T 0(\I 1) = \I 0$.

\medskip \noindent (ii) $\imply $ (iii): Obvious.

\medskip \noindent (iii) $\imply $ (i): Follows immediately from \ref
{InclusionInduced}.

\endProof

In order to give a concrete description of a transposed ideal, we must bring in
certain important maps between the various group algebras in sight.  Initially,
consider the natural projection
  $$
  P: \textstyle \sum _{g\in G}c_g\delta _g \in \GpAlg G \mapsto \sum _{h\in \H
0}c_h\delta _h \in \Ga 0,
  $$
  and, given $k$ and $l$ in $G$, define the map
  $$
  \pconj _{k,l}: c \in \Ga 1 \mapsto P(\delta _{k\inv }c\delta _l)\in \Ga 0.
  $$
  For an explicit expression, let $c=\sum _{h\in \H 1}c_h\delta _h\in \Ga 1$,
and notice that
  $$
  \pconj _{k,l}(c) =
  P\Big (\sum _{h\in \H 1}c_h\delta _{k\inv hl}\Big ) =
  \sum \subsub {h\in \H 1}{k\inv hl\in \H 0}c_h\delta _{k\inv hl} =
  \delta _{k\inv }\Big (\sum _{h\in \H 1\cap \,k\H 0l\inv }c_h\delta _h\Big
)\delta _l.
  $$

\state Proposition \label DescribeTransposition
  Let $\I 0$ be an admissible ideal in $\Ga 0$.  Then the transposition of $\I
0$ to $\Ga 1$ is given by
  $$
  \T 1(\I 0) = \bigcup _{V\ni \x 1}\bigcap \subsub {k,l\in \S 0}{\theta _k(\x
0)\in V}\pconj _{k,l}\inv (\I 0),
  $$
  where by ``$\,V\ni \x 1$'' we mean that $V$ ranges in the family of all
neighborhoods of $\x 1$.

\Proof
  Let $c=\sum _{h\in \H 1}c_h\delta _h \in \Ga 1$.  Then by \ref {LiftC} one has
that $c$ lies in $\T 1(\I 0)= \F 1\big (\J 0{\I 0}\big )$ if and only if there
exists a {\klopen } set $V$, such that
  $$
  \x 1\in V\subseteq X_h,
  \equationmark CondOnV
  $$
  whenever $c_h\neq 0$, and
  $$
  c_V:= \sum _{h\in \H 1}c_h1_V\d _h \in \J 0{\I 0}.
  $$

Using \ref {UseDeltaKbDeltaL}, the above is equivalent to saying that, for every
$k$ and $l$ in $\S 0$, one has that
  $$\def \quad {\ \ }\def \crr {\stake {19pt}\hfill \ = \cr }
  \matrix {\I 0 & \ni & \ds \sum _{h\in \H 1\cap \,k\H 0 l\inv }c_h1_V\big
(\theta _k(\x 0)\big )\, \delta _{k\inv hl} \crr
              & = & 1_V\big (\theta _k(\x 0)\big ) \kern -5pt \ds \sum _{h\in \H
1\cap \,k\H 0 l\inv }c_h\, \delta _{k\inv hl} \crr
              & = & 1_V\big (\theta _k(\x 0)\big ) \ \pconj _{k,l}(c).\hfill }
  $$
  This condition is clearly meaningless unless $\theta _k(\x 0)$ is in $V$, in
which case it says that $c\in \pconj _{k,l}\inv (\I 0)$.

  If follows that $c\in \T 1(\I 0)$ if and only if $c$ lies in the set whose
definition is almost exactly what the statement claims $\T 1(\I 0)$ to be, the
only difference being the family of sets where $V$ ranges which, in the present
case, consists of all {\klopen } neighborhoods of $\x 1$ satisfying \ref
{CondOnV}.  However, if we take into account that $\x 1$ admits a fundamental
system of {\klopen } neighborhoods, and that the correspondence
  $$
  V\mapsto \bigcap \subsub {k,l\in \S 0}{\theta _k(\x 0)\in V}\pconj _{k,l}\inv
(\I 0),
  \equationmark OneVatATime
  $$
  is decreasing, then we see that such a difference is irrelevant.  \endProof

An interesting consequence is that, when $\x 1$ is not in the closure of the
orbit of $\x 0$, the transposition of ideals leads to a triviality:

\state Proposition Suppose that $\x 1\notin \cOrb (\x 0)$.  Then, for every
ideal $\I 0\ideal \Ga 0$, one has that $\T 1(\I 0)=\Ga 1$.

\Proof Let $V$ be a neighborhood of $\x 1$ such that $V\cap \cOrb (\x
0)=\emptyset $. Then there is no $k$ in $\S 0$ such that $\theta _k(\x 0)\in V$,
whence \ref {OneVatATime} consists of the intersection of the empty family of
sets, resulting in the universe where it is considered, namely $\Ga 1$.
\endProof

The transposition towards strongly regular points may be described in a much
simpler way:

\state Theorem
  Assume that $\x 1$ is strongly regular, and let $\I 0$ be an admissible ideal
in $\Ga 0$.  Then
  $$
  \T 1(\I 0) =
  \bigcup _{V\ni \x 1}\bigcap \subsub {k\in \S 0}{\theta _k(\x 0)\in V}\delta
_k\I 0 \delta _{k\inv }.
  $$

\Proof
  Let
  $
  c=\sum _{h\in \Gamma }c_h\delta _h \in \Ga 1,
  $
  where $\Gamma $ is a finite subset of $\H 1$.  Then by \ref {LiftC} we have
that the following two conditions are equivalent:
  \izitem
  \zitem $c\in \T 1(\I 0)= \F 1\big (\J 0{\I 0}\big )$,
  \zitem there exists a {\klopen } set $V$, such that
  $$
  \x 1\in V\subseteq X_h, \for h\in \Gamma ,
  $$
  and
  $$
  c_V:= \sum _{h\in \Gamma }c_h1_V\d _h \in \J 0{\I 0}.
  $$

  Let us prove that (ii) is in turn equivalent to:

\bigskip

\zitem there exists a {\klopen } set $V$, satisfying all of the requirements of
(ii), and morever such that $\theta _h$ fixes $V$, for every $h$ in $\Gamma $.

\bigskip To see that (ii) implies (iii), use the fact that $\x 1$ is strongly
regular to produce a {\klopen } neighborhood $W$ of $\x 1$, such that $\theta
_h$ fixes $W$, for every $h$ in $\Gamma $.  One then has that
  $$
  \J 0{\I 0} \ni 1_Wc_V = \sum _{h\in \Gamma }c_h1_W1_V\d _h = \sum _{h\in
\Gamma }c_h1_{W\cap V}\d _h = c_{W\cap V},
  $$
  thus proving (iii).  That (iii) implies (ii) is evident.

Assuming that $c\in \T 1(\I 0)$, and hence that (iii) holds, it follows from
\ref {LemaOnMembership} that, for every $k$ in $\Sx $, with $\theta _k(x_0)\in
V$, one has that
  $
  \delta _{k\inv }c\delta _k\in I.
  $
  Consequently $c\in \delta _kI\,\delta _{k\inv }$, which is to say that
  $$
  c\in \bigcap \subsub {k\in \S 0}{\theta _k(\x 0)\in V}\delta _k\I 0 \delta
_{k\inv },
  $$
  which in turn implies that $c$ belongs to the set the statement claims $\T
1(\I 0)$ to be.

Conversely, if $c$ lies in that set, there exists an open neighborhood $V$ of
$\x 1$ such that, whenever $k\in \S 0$, and $\theta _k(\x 0)\in V$, one has that
$c\in \delta _kI\,\delta _{k\inv }$.  Since $\x 1$ is strongly regular, and upon
shrinking $V$ if necessary, we may suppose that $V$ is {\klopen }, and that
$\theta _h$ fixes $V$, for every $h$ in $\Gamma $.  It then follows from \ref
{LemaOnMembership} that $c_V\in \J 0{\I 0}$, namely that condition (ii) above
holds, so that (i) also holds, so $c\in \T 1(\I 0)$.  This completes the proof.
\endProof

\references

\Article LisaEtAll
  J. Brown, L. Clark, C. Farthing and A, Sims;
  Simplicity of algebras associated to \'etale groupoids;
  Semigroup Forum, 88 (2014), 433-452

\Bibitem ReznikofEtAll
  J. H. Brown, G. Nagy, S. Reznikoff and A. Sims;
  Cartan subalgebras in C*-algebras of Hausdorff etale groupoids;
  preprint, arXiv:1503.03521v3 [math.OA], 2016

\Bibitem EH
  E. G. Effros, F. Hahn;
  Locally compact transformation groups and C*-algebras;
  Memoirs of the American Mathematical Society, no. 75, 1967

\Bibitem mybook
  R. Exel;
  Partial Dynamical Systems, Fell Bundles and Applications;
  to be published in a forthcoming NYJM book series.  Available from {\tt
  http://mtm.ufsc.br/$\sim $exel/papers/pdynsysfellbun.pdf}

  \def \Bibitem #1 #2; #3; #4 \par {\smallbreak
    \global \advance \bibno by 1
    \item {[\possundef {#1}]} #2, {``#3''}, #4.\par
    \ifundef {#1PrimarilyDefined}\else
      \fatal {Duplicate definition for bibliography item ``{\tt #1}'', already
defined in ``{\tt [\csname #1\endcsname ]}''.}
      \fi
	\ifundef {#1}\else
	  \edef \prevNum {\csname #1\endcsname }
	  \ifnum \bibno =\prevNum \else
		\error {Mismatch bibliography item ``{\tt #1}'', defined earlier (in aux file
?) as ``{\tt \prevNum }'' but should be
	``{\tt \number \bibno }''.  Running again should fix this.}
		\fi
	  \fi
    \define {#1PrimarilyDefined}{YES}\if \TRUE \auxfile \immediate \write 1
{\textbackslash newbib {#1}{\number \bibno }}\fi }

\Article Dan
  D. Gon\c {c}alves, J. \umlaut Oinert and D. Royer;
  Simplicity of partial skew group rings with applications to Leavitt path
algebras and topological dynamics;
  J. Algebra, 420 (2014), 201-216

\Article GR
  E. C. Gootman and J. Rosenberg;
  The structure of crossed product C*-algebras: a proof of the generalized
Effros-Hahn conjecture;
  Invent. Math., 52 (1979), no. 3, 283-298

\Article IonWill
  M. Ionescu and D. Williams;
  The generalized Effros-Hahn conjecture for groupoids;
  Indiana Univ. Math. J., 58 (2009), no. 6, 2489-2508

\Article Renault
  J. Renault;
  The ideal structure of groupoid crossed product C*-algebras. With an appendix
by Georges Skandalis;
  J. Operator Theory, 25 (1991), no. 1, 3-36

\Article Sauvageot
  Jean-Luc Sauvageot;
  Ideaux primitifs de certains produits croises;
  Math. Ann., 231 (1977), 61-76

\Article FS
  T. Fack, G. Skandalis;
  Sur les repr\'esentations et id\'eaux de la C*-alg\`ebre d'un feuilletage;
  J. Oper. Theory, 8 (1983), 95-129

\Bibitem SimsWill
  A. Sims and D. Williams;
  The primitive ideals of some \'etale groupoid C*-algebras;
  preprint, arXiv:1501.02302 [math.OA], 2015

\Article SteinbAlg
  B. Steinberg;
  A groupoid approach to discrete inverse semigroup algebras;
  Adv. Math., 223 (2) (2010), 689-727

\Article Steinberg
  B. Steinberg;
  Simplicity, primitivity and semiprimitivity of \'etale groupoid algebras with
applications to inverse semigroup algebras;
  J. Pure Appl. Algebra, 220 (2016), 1035-1054

\Bibitem Willard
  S. Willard;
  General topology;
  Addison-Wesley, 1970, xii+369 pp

\Bibitem WillBook
  D. Williams;
  Crossed products of C*-algebras;
  Mathematical Surveys and Monographs, vol. 134, American Mathematical Society,
2007

\endgroup  \bye